\documentclass[12 pt]{article}
\usepackage{amssymb,amsmath,amsfonts,amsthm}
\newcommand\BBP{{\mathbb {P}}}
\newcommand\BBI{{\mathbb {I}}}

\newcommand\R{{\mathbb {R}}}
\newcommand\Z{{\mathbb {Z}}}
\newcommand\F{{\mathcal {F}}}

\newcommand\BBE{{\mathbb {E}}}

\newtheorem {Lemma}{Lemma}[section]
\newtheorem {Theorem}{Theorem}[section]
\newtheorem {Proposition}{Proposition}[section]
\newtheorem {Corollary}{Corollary}[section]
\newtheorem {Fact}{Fact}[section]

\theoremstyle{definition}
\newtheorem{Definition}{Definition}[section]
\newtheorem{Notation}{Notation}[section]
\newtheorem{Remark}{Remark}[section]
\newtheorem{Comment}{Comment}[section]

\newcommand\I{{ 1\hspace{-1,2mm}{\mathrm I}}}

\newcommand\no{\noindent}

\newcommand\beq{\begin{equation}}
\newcommand\eeq{\end{equation}}

\def\S{\tilde S}
\def\X{\tilde X}
\def\no{\noindent}

\def\cov{\mathop{\rm Cov}\limits}

\def\XN{(X_i)_{i\in{\mathbb N}}}
\def\XZ{(X_i)_{i\in{\mathbb Z}}}

\begin{document}
\title{Quadratic transportation cost in the conditional central limit theorem for dependent sequences}

\author{J\'er\^ome Dedecker\footnote{J\'er\^ome Dedecker, Universit\'e de Paris, CNRS, MAP5, UMR 8145,
45 rue des  Saints-P\`eres,
F-75006 Paris, France.},
Florence Merlev\`ede\footnote{Florence Merlev\`ede, LAMA,  Univ Gustave Eiffel, Univ Paris Est Cr\'eteil, UMR 8050 CNRS,  \  F-77454 Marne-La-Vall\'ee, France.} 
and Emmanuel Rio \footnote{Emmanuel Rio, Universit\'e de Versailles, LMV UMR 8100 CNRS,  \ 45 avenue des Etats-Unis, 
F-78035 Versailles, France.}
}

\maketitle
\begin{abstract}
In this paper, we give estimates of  the quadratic transportation cost in the conditional central limit theorem for a large class of dependent sequences. Applications to irreducible Markov chains, dynamical systems generated by intermittent maps and $\tau$-mixing sequences are given.  \\

\noindent{\it  MSC2020 subject classifications: }   60F05; 60F25; 60E15; 37E05. 
\no{\it Running head: }Quadratic transportation cost in the conditional CLT
\end{abstract}

\section{Introduction} \label{Intro}
\setcounter{equation}{0}

Let $\XZ$ be a strictly
stationary sequence of real-valued random variables (r.v.) with mean
zero and finite variance. Set $S_n = X_1 + X_2 + \cdots + X_n$. By
$P_{n^{-1/2}S_n}$ we denote the law of $n^{-1/2}S_{n}$ and by
$G_{\sigma^2}$  the normal distribution $ N(0, \sigma^2)$. In this
paper, we assume furthermore that the series $\sigma^2=\sum_{k \in {\mathbb Z}} {\rm Cov} (X_0,X_k)$ is convergent (under this assumption $\lim_n n^{-1} {\rm Var} S_n = \sigma^2$) and we shall give  quantitative estimates of the approximation of $P_{n^{-1/2}S_n}$ by $G_{\sigma^2}$ in terms of the quadratic cost,
which is the square of the  ${\mathbb L}^2$-minimal distance.  With this aim, we first recall the definition of the 
${\mathbb L}^p$-minimal metrics. 

Let ${\mathcal L}(\mu, \nu)$ be the set of  probability laws on
$\mathbb R^2$ with  marginals $\mu$ and $\nu$. For $p\geq 1$, let
\[
W_p(\mu, \nu) =  \displaystyle \inf \Big \{ \Big (\int |x-y|^p P(dx, dy) \Big )^{1/p}
: P \in {\mathcal L}(\mu, \nu) \Big \}   \, .\]
$W_p$ is usually called the ${\mathbb L}^p$-minimal distance, and sometimes  the Wasserstein 
distance of order $p$.  It is well known that for probability laws $\mu$ and $\nu$ on $\mathbb R$ 
with respective distributions functions (d.f.) $F$ and $G$,
\begin{equation} \label{def2wasser}
W_p(\mu, \nu) =  \Big ( \int_0^1 |F^{-1}(u) - G^{-1}(u) |^p du\Big
)^{1/p} \, ,
\end{equation}
where $F^{-1}$ and $G^{-1}$ denote respectively the generalized inverse functions of $F$ and $G$.  We refer to Chapter 6 in Villani \cite{Vi2009} for the properties of this metric. 

For  $\XZ$  a sequence of independent and identically distributed (iid) centered real valued random variables  in ${\mathbb L}^4$, with  variance $\sigma^2$,  inequality (1.7) in Rio \cite{Rio09} states that there exists a universal constant $c$ such that for any positive integer $n$
\beq \label{ineRioUB}
n W_2^2 (P_{S_n / \sqrt n} ,  G_{\sigma^2}\big )  \leq  c \,  \sigma^{-2}  \Vert X_1 \Vert_4^4 \, .
\eeq
In addition, it is also shown in the same paper that this upper bound is optimal. More precisely, for any 
$\kappa \geq 1$, let ${\cal M}(4, \kappa)$ be the class of the probability measures  $\mu$ on the real line such that $\int  x d \mu (x) = 0$, $\int  x^2 d \mu (x) = 1$ and  $\int  x^4 d \mu (x) = \kappa$.  In case where $\XZ$ is a sequence of iid random variables with common law $\mu $ in  ${\cal M}(4, \kappa)$,  Theorem  5.1 in \cite{Rio09} asserts that 
\beq \label{ineRioLB}
 \sup_{ \mu \in  {\cal M}(4, \kappa)}  \liminf_{n \rightarrow \infty} n W_2^2 (P_{S_n / {\sqrt n} } ,  G_{1}\big )  \geq  \kappa /12 \, .
\eeq

We refer to Bobkov \cite{Bob13} for another proof of  \eqref{ineRioUB} based on relative entropy and Talagrand's entropy-transport inequality.  Actually, the following more general  result   holds: for any $p \geq 1$, there exists a universal constant $c_p$ such that for any positive integer $n$, 
\[
n^{p/2} W_p^p (P_{S_n / \sqrt n} ,  G_{\sigma^2}\big )  \leq  c_p \,  \sigma^{-p}  \Vert X_1 \Vert_{p+2}^{p+2} \, .
\]
(see Rio \cite{Rio09} for $p \in [1,2]$ and Bobkov \cite{Bob18} for $p>2$).  Extensions to random vectors in ${\mathbb R}^d$ are given in Bonis \cite{Bo20}.  We also mention the extensions of  the upper bound \eqref{ineRioUB} to the $m$-dependent case and to $U$-statistics  obtained by  Fang \cite{Fa19}.

In this paper, one of our motivations is to relax the independence assumption and to find sufficient conditions in case of dependent sequences ensuring that 
\beq  \label{conj1}
W_2 (P_{S_n / \sqrt n } ,  G_{\sigma^2}\big )  = O (n^{-1/2}) \, .
\eeq   
In the dependent setting, a well known class is the class of irreducible  aperiodic  and positively recurrent Markov chains $(\xi_n)$ with an atom denoted by $A$ 
(see the definition page 286 in \cite{Bo82}).  Let $\pi $ be the unique  invariant distribution of the Markov chain. From now on, $(\xi_n)$ will be the Markov chain starting from  $\pi$. Let us then consider the strictly stationary sequence $(X_k) $ defined by $X_k =f ( \xi_k)$ with $f$ a bounded function such that $\pi(f) =0$. In view of  the regeneration scheme and the upper bound \eqref{ineRioUB},  one can conjecture that \eqref{conj1}  holds for $S_n = \sum_{k=1}^n X_k$ provided that $\BBE_A ( \tau_A^4) < \infty$ where $\tau_A$ is the first return time in $A$ and ${\BBE}_A$ stands for the expectation under ${\BBP}_x$ for $x \in A$.  Next, from   \cite[Lemma 3]{Bo82} and \cite[page 165]{Rio17}, it is known  that  $\BBE_A ( \tau_A^4) < \infty$  is equivalent to 
\beq \label{condalpha4}
 \sum_{n>0} n^{2} \alpha_n < \infty \, , 
\eeq
where $
\alpha_n = { \textstyle  \frac{1}{4}  } \sup_{\Vert f \Vert_{\infty} \leq 1}    \Vert  \BBE ( f(\xi_n ) | \xi_0 ) - \BBE (f (\xi_n) )  \Vert_1  
$.

In this paper we shall prove that \eqref{conj1} holds true for any stationary sequence $(X_k)_{k  \in {\mathbb Z}}$  of  bounded real-valued random variables  satisfying  \eqref{condalpha4} for the sequence  $(\alpha_n)_{n \geq 0}$ of strong mixing coefficients in the sense of Rosenblatt (see for instance \cite[Section 5.1.1.]{MPU} for a definition of these coefficients in the general case), which includes the case of Markov chains described above.  This will be a consequence of a more general result also valid for a class of weakly dependent sequences,   which may fail to be strongly mixing. In order to give more precise statements of our results, let us now introduce the  dependence coefficients that we will use in this paper. 

\begin{Definition} \label{deftheta} Let $\XZ$ be a stationary sequence of bounded real-valued random variables and ${\mathcal  F}_0 = \sigma ( X_i, i \leq 0)$. Let 
$ \Gamma_{p,q} = \{  (a_i)_{ 1 \leq i \leq p}   \in {\mathbb N}^p \, : \,  a_1 \geq 1 \text{ and } \sum_{i=1}^p  a_i \leq q  \} $, for $p$ and $q$ positive integers. For $k\geq 0$, set
\[
\theta_{X,p,q} (k)  = \sup_{ k_p>k_{p-1}> \ldots > k_2>k_1 \geq k \atop (a_1, \dots, a_p ) \in \Gamma_{p,q}} \Big \Vert  \BBE \Big ( \prod_{i=1}^p X_{k_i}^{a_i} |   {{\mathcal  F}_0 }  \Big  ) - \BBE  \Big ( \prod_{i=1}^p X_{k_i}^{a_i}  \Big )  \Big \Vert_1 \, .
\]
\end{Definition}
As a consequence of our Theorem \ref{thW2},  we will obtain that if 
\beq  \label{condtheta4}
\sum_{k \geq 1} k^2 \theta_{X,4,4} (k) < \infty \, , 
\eeq 
 then \eqref{conj1}  holds, which immediately implies that  \eqref{conj1}  holds for additive bounded functionals of a Markov chain  satisfying  \eqref{condalpha4}. In fact we shall give a conditional version of \eqref{conj1}  and show that when $(X_k)_{k  \in {\mathbb Z}}$ is a stationary sequence of centered and bounded real-valued random variables satisfying \eqref{condtheta4} then 
\beq  \label{condW2res}
 \BBE \big ( W_2^2 (P_{S_n / {\sqrt n} |{\mathcal F}_0 } ,  G_{\sigma^2}) \big ) = O ( n^{-1}  ) \, . 
\eeq
Note that in case of bounded functions of a Markov chain $(\xi_k)_k$  satisfying $\BBE_A ( \tau^4) < \infty$, with invariant distribution $\pi$, the Schwarz inequality together with  \eqref{condW2res}  imply that 
\[ 
\BBE_\mu \big ( W_2 (P_{S_n / {\sqrt n} |{\xi}_0 } ,  G_{\sigma^2}) \big ) = O ( n^{-1/2}  )
\]
for any positive measure $\mu$ such that $d \mu = f d \pi $ with $\int f^2 d \pi < \infty$.  Above $ \BBE_\mu$ stands for the expectation of the chain under the initial law $\mu$. 

\smallskip

It is noteworthy to indicate that \eqref{condW2res} implies \eqref{conj1}. Indeed the following fact is valid. 

\begin{Fact} 
\label{CompQCCondQC}
Let  $X$ and $Y$ be two random variables defined on $(\Omega,{\mathcal A}, {\mathbb P})$ and ${\mathcal F}$ be a sub $\sigma$-algebra of ${\mathcal A}$. Then
$W_2^2  (P_{X }  ,  P_{Y }  ) \leq \BBE \big (  W_2^2  (P_{X | {\mathcal F}}  ,  P_{Y | {\mathcal F}}  )   \big )$.
\end{Fact}

To see this, let $U$ be a random variable with uniform distribution over $[0,1]$, independent of ${\mathcal F}$,  and let 
$F_{ X| {\mathcal F} }$ and $F_{ Y| {\mathcal F}  }$ denote respectively the conditional distribution functions of 
$X$ and $Y$ given ${\mathcal F}$. 
Set $X^* = F_{ X| {\mathcal F} }^{-1} (U)$ and $Y^* = F_{Y| {\mathcal F} }^{-1} (U)$.  
Then $X^*$ has the law $P_X$,  $Y^*$ has the law $P_Y$ and, by   \eqref{def2wasser}, 
$W_2^2  (P_{X | {\mathcal F}}  ,  P_{Y | {\mathcal F}}  )  = \BBE \big ( |X^* - Y^*|^2 | {\mathcal F}  \big )$.
Taking the expectation, it implies the above fact, since $W_2^2$ is the minimal quadratic cost.
 \par\medskip
 
 To prove Theorem 2.1, we shall apply  Lindeberg's method, which was used by Billingsley \cite{Bil61} and Ibragimov
 \cite{Ibr63} in the case of martingales with stationary differences to prove the central limit theorem (we also consider this particular case in our Theorem \ref{thW2M}). Note that this method was adapted to a large class of dependent sequences (non necessarily martingale differences) to evaluate the  ${\mathbb L}^1$-minimal  distance between $P_{S_n / {\sqrt n}}$ and $G_{\sigma^2}$,  by P\`ene \cite{Pe05} in the bounded multidimensional case, and next  by Dedecker and Rio \cite{DR08} in the unbounded case (under conditions involving  some coefficients similar to $\theta_{X,4,3}$, or weak mixing coefficients such as those described in Definition \ref{defalpha} below). Recently, estimates of the ${\mathbb L}^1$-minimal  distance between $P_{S_n / {\sqrt n}}$ and $G_{\sigma^2}$ when the underlying process is a function of  iid random variables are given in Theorem 3.1 in \cite{JWZ}. Their  conditions are expressed in terms of some  coupling coefficients.

Our paper is organized as follows. Section \ref{section2} is devoted to the statements of upper bounds concerning the quadratic transportation cost in the conditional central limit theorem and their applications to pointwise estimates
for the distribution function of the normalized sums and its generalized inverse. 
Applications to  $\alpha$-dependent sequences, $\tau$-mixing sequences and symmetric random walk in the circle  are given in Section \ref{section3}. The proofs are postponed to Section \ref{section4}. 
Links between $ | F_{S_n/\sigma_n}^{-1}(u)  -  \Phi^{-1} (u) | $ and $W_2 (P_{S_n / \sigma_n} ,  G_1)$ are given in Section \ref{Annex}, where  $\sigma_n = \sqrt{ {\rm Var} S_n}$, 
$\Phi^{-1}$ is the  inverse  of the distribution function of the  standard normal distribution and 
$F_{S_n/\sigma_n}^{-1}$ is the generalized inverse of the distribution function of $S_n/\sigma_n$. 

\medskip

In the rest of the paper, we shall use the following notation:  for two sequences $(a_n)_{n \geq 1}$ and $(b_n)_{n \geq 1}$ of positive reals, $a_n \ll b_n$ means there exists a positive constant $C$ not depending on $n$ such that $a_n \leq C b_n$ for any $n\geq 1$.  Moreover, for a real-valued  random variable $X$ in ${\mathbb L}^1$, the notation $X^{(0)}$ means 
$X - \BBE (X)$. 

\section{Quadratic cost in the conditional CLT} \label{section2}
\setcounter{equation}{0}

The main result of this paper is Theorem \ref{thW2} below. 

\begin{Theorem} \label{thW2} Assume that $\Vert X_0 \Vert_{\infty} \leq M$ and that $\sum_{k \geq 1}  \theta_{X,2,2} (k) < \infty$. Then $\sigma^2 = \BBE (X_0^2) + 2  \sum_{k \geq 1}\BBE (X_0X_k)$ converges and 
$$
\BBE \big ( W_2^2 (P_{S_n / {\sqrt n} |{\mathcal F}_0 } ,  G_{\sigma^2}) \big ) \ll   
n^{-1/2} \Big ( 1 + \sum_{k\geq 1}  ( k \wedge \sqrt{n} )  \theta_{X,2,2} (k) \Big ) \, .\eqno (a)
$$
If furthermore  $ \sum_{k\geq 1} k  \theta_{X,4,4} (k) < \infty$, then
$$
\BBE \big ( W_2^2 (P_{S_n / {\sqrt n} |{\mathcal F}_0 } ,  G_{\sigma^2}) \big ) \ll   n^{-1} \Big ( 1 + \sum_{k\geq 1} k ( k \wedge \sqrt{n} )  \theta_{X,4,4} (k) \Big )\, .\eqno (b)
$$
\end{Theorem}

\begin{Comment} Item (a) provides a rate in the CLT for the $W_2$-metric as soon as $\sum_{k \geq 1}  \theta_{X,2,2} (k) < \infty$. In addition, if $ \sum_{k\geq 1}  k \theta_{X,2,2} (k) < \infty$, then the rate  in the $W_2$-metric is of order $n^{-1/4}$.  Furthermore, by Item (b), if $ \sum_{k\geq 1} k  \theta_{X,4,4} (k) < \infty$, then  the rate  in the CLT for the $W_2$-metric is  $o ( n^{-1/4})$. For example, if $\theta_{X,4,4} (k)  = O (k^{-a}) $ with $ a \in ]1,3[ $ and $a \neq 2$,   Theorem \ref{thW2} implies that $W_2 (P_{S_n / {\sqrt n} } ,  G_{\sigma^2}) \ll     n^{-(a-1)/4} $.  Moreover $  W_2 (P_{S_n / {\sqrt n} } ,  G_{\sigma^2}) \ll  n^{-1/2}$ as soon as $ \sum_{k\geq 1}  k^2  \theta_{X,4,4} (k) < \infty$.
\end{Comment}

\begin{Comment} \label{commentaftermainth} Assume $\sigma >0$. Set $\sigma_n= \sqrt{ {\rm Var} S_n }$. If $\sigma>0$, then $\sigma_n>0$ for any positive $n$. Set  
$\kappa_2 = \BBE \big (W_2^2 (  P_{S_n / \sigma {\sqrt n}  |{\mathcal F}_0 } ,  P_{S_n / \sigma_n  |{\mathcal F}_0 } )\big )$:
\[
\kappa_2    =  
\Big ( \frac{\sigma_n}{\sigma \sqrt n}  - 1 \Big )^2  \leq  \Big ( \frac{\sigma_n^2}{n\sigma^2}  - 1 \Big )^2  
 = 
(n\sigma^2)^{-1}  \Big| \frac{\sigma_n^2}{n\sigma^2}  - 1 \Big| \,  \Big| 2 \sum_{k \geq 1} (k \wedge n)  {\rm Cov} (X_0,X_k)\Big| . 
\]
Now, from the definition of the coefficients $ \theta_{X,1,1} (k)$,
\beq \label{UpperBoundCov1}
\sum_{k \geq 1} (k \wedge n) | {\rm Cov} (X_0,X_k) | 
\leq \Vert X_0 \Vert_{\infty} \sum_{k \geq 1} (k \wedge n ) \theta_{X,1, 1} (k) . 
\eeq
Therefore, if in addition $\Vert X_0 \Vert_{\infty} \leq M$, 
$\kappa_2    \ll n^{-1} M \Vert_{\infty} \sum_{k \geq 1} (k \wedge n ) \theta_{X,1, 1} (k)$, 
which is always of a smaller order than the upper bounds (a) and (b).  Hence Theorem \ref{thW2}  also holds for 
$\BBE \big ( W_2^2 (P_{S_n / \sigma_n |{\mathcal F}_0 } ,  G_{1}) \big ) $. 
\end{Comment}

We now give applications of Theorem \ref{thW2} to pointwise estimates. We start by Berry-Esseen type estimates. Arguing for instance  as in \cite[Remark 2.4]{DMR09}, Theorem \ref{thW2} together with Comment 2.2 
imply the following upper bound.

\begin{Corollary} \label{AppliBE}  Assume that $\sigma>0$, $\Vert X_0 \Vert_{\infty} \leq M$ and $\sum_{k \geq 1} k^2  \theta_{X,4,4} (k) < \infty$. Then 
\[   \Delta_{n} =  \sup_{x \in {\mathbb R}}\big |  \BBP ( S_n/\sigma_n  \leq x ) - \Phi( x )\big |\ll   n^{-1/3}  \, .
\]
\end{Corollary}

We now give applications of our main result to estimates of the quantiles and the superquantiles
of $S_n/\sigma_n$ in the nondegenerate case.   Define the $1$-risk $Q_{1,X}$ of $X$, as in Pinelis \cite{Pi14},  by 
\beq
\label{DefSuperQuantile}
Q_{1,X} (u) = \frac 1 u  \int_0^u F_X^{-1} (1-t) dt .
\eeq
Then $Q_{1,X} (u)$ is the value of the  superquantile of $X$ at point $(1-u)$. The corollary below, which is a  
consequence of Theorem \ref{thW2} and Proposition \ref{PropVaR} provides estimates 
of the accuracy in the central limit theorem for $F^{-1}_{S_n/\sigma_n}$ and $Q_{1,S_n/\sigma_n}$.  Its proof is given in Section \ref{Annex}. 

\begin{Corollary} \label{AppliVaR}  Assume that $\Vert X_0 \Vert_{\infty} \leq M$,  $\sum_{k \geq 1} k^2  
\theta_{X,4,4} (k) < \infty$ and $\sigma^2 >0$. Let $Y$ be a standard normal. 
Then there exists some constant $C>0$ such that, 
for any  $n\geq 1$ and any $u$ in $(0,1)$,  
$$
| F_{S_n/\sigma_n}^{-1} (u) - \Phi^{-1} (u) | \leq 
C \max \bigl( (nu(1-u))^{-1/2} , ( nu(1-u) )^{-1/3} |\log (u(1-u) ) |^{-1/6} \bigr)  \eqno (a)
$$
and 
$$
| Q_{1, S_n/\sigma_n} (u) - Q_{1,Y} (u) | \leq C (nu)^{-1/2} \sqrt{1-u} . 
\eqno (b) 
$$
\end{Corollary}

\begin{Comment} From Corollary \ref{AppliVaR}(a), for any sequence $(\varepsilon_n)_n$ of reals in $(0,1/2)$ 
such that $\lim_n \varepsilon_n= 0$ and $\lim_n n\varepsilon_n = \infty$,
$$
\lim_{n\rightarrow \infty} 
\sup_{u \in [\varepsilon_n ,1-\varepsilon_n] } | F_{S_n/\sigma_n}^{-1} (u) - \Phi^{-1} (u) | = 0 ,
$$ which can not be deduced from a Berry-Esseen type bound with the rate $n^{-1/2}$.  Indeed, 
if $\Delta_n$ is defined as  in Corollary \ref{AppliBE}, one can only get that
\[
| F_{S_n/\sigma_n}^{-1} (u) - \Phi^{-1} (u) |   \leq      \Phi^{-1} (\min (1, u + \Delta_n) ) - \Phi^{-1} (u)  
\]
for $u\geq 1/2$, which is of interest only if  $u < 1 - \Delta_n$. 
%Similarly, Corollary \ref{AppliVaR}(b) implies that 
%$$
%\lim_{n \rightarrow \infty} \sup_{u \in [\varepsilon_n ,1] } | Q_{1,S_n/\sigma_n} (u) - Q_{1,Y} (u) | = 0 ,
%$$
% which cannot be deduced from  the rate $n^{-1/2}$ in the CLT for the $W_1$-metric or a $W_p$-metric for $p<2$. 
\end{Comment}
 
 If furthermore the sequence $\XZ$  is a sequence of martingale differences, then the conditions on the 
 dependence coefficients can be weakened as follows (the proof being less intricate  is left to the reader).  
 
\begin{Theorem} \label{thW2M} Assume that $\XZ$ is a sequence of martingale differences such that $\Vert X_0 \Vert_{\infty} \leq M$ and $\BBE (X_0^2) = \sigma^2$. Then
$$
\BBE \big ( W_2^2 (P_{S_n / {\sqrt n} |{\mathcal F}_0 } ,  G_{\sigma^2}) \big ) \ll   n^{-1/2} 
\Big ( 1 + \sum_{k=1}^ {[\sqrt n]}   \theta_{X,1,2} (k) \Big ) \, .  \eqno (a)
$$
If furthermore  $ \sum_{k\geq 1}  \theta_{X,3,4} (k) < \infty$, then
$$
\BBE \big ( W_2^2 (P_{S_n / {\sqrt n} |{\mathcal F}_0 } ,  G_{\sigma^2}) \big ) \ll   n^{-1} \Big ( 1 + \sum_{k\geq 1}  ( k \wedge \sqrt{n} )  \theta_{X,3,4} (k) \Big )\, .  \eqno (b) $$
%In particular,
%\begin{itemize}
%\item[1.] If $\theta_{X,4} (k)  = O (k^{-a}) $ with $2 < a < 3$,  $\BBE \big ( W_2^2 (P_{S_n / {\sqrt n} |{\mathcal F}_0 } ,  G_{\sigma^2}) \big ) \ll     n^{-(a-1)/2} $. \\
%\item[2.] If $\sum_{k \geq 1} k^2 \theta_{X,4} (k) < \infty$, $\BBE \big ( W_2^2 (P_{S_n / {\sqrt n} |{\mathcal F}_0 } ,  G_{\sigma^2}) \big ) \ll    n^{-1} $. \end{itemize} 
\end{Theorem}

\begin{Comment} Item (a) provides a rate in the CLT as soon as $\theta_{X,1,2} (k) = o (1)$. If $\theta_{X,1,2} (k)  = O (k^{-a}) $ with $a$ in $(0,1) $, (a) ensures that $W_2 (P_{S_n / {\sqrt n} } ,  G_{\sigma^2}) \ll     n^{-a/4} $. If $ \sum_{k\geq 1}  \theta_{X,1,2} (k) < \infty$, then the rate  is of order $n^{-1/4}$. 
Item (b) provides faster rates  under the condition $ \sum_{k\geq 1}  \theta_{X,3,4} (k) < \infty$. 
Indeed the rate of convergence under this condition is  $o ( n^{-1/4})$. 
If $\theta_{X,3,4} (k)  = O (k^{-a}) $ with $a$ in $(1,2)$, (b) ensures that $W_2 (P_{S_n / {\sqrt n} } ,  G_{\sigma^2}) \ll     n^{-a/4} $.  Moreover $  W_2 (P_{S_n / {\sqrt n} } ,  G_{\sigma^2}) \ll  n^{-1/2}$ as soon as 
$ \sum_{k\geq 1}  k  \theta_{X,3,4} (k) < \infty$.
\end{Comment}

\section{Examples} \label{section3}

\setcounter{equation}{0}

\subsection{$\alpha$-mixing sequences}\label{alphamixing}
Let $(\Omega, {\mathcal A}, \BBP)$ be a probability space and let $ {\mathcal U}$ and $ {\mathcal V}$
 be two $\sigma$-algebras of $ {\mathcal A}$.  The strong mixing coefficient $ \alpha (\mathcal{U},\mathcal{V})  $ between these $\sigma$-algebras is defined as follows: 
 \[
  \alpha (\mathcal{U},\mathcal{V})  =\sup \{ | {\ \mathbb{P}} (U\cap
V) - {\mathbb{P}} (U) { \mathbb{P}} (V)|: U \in \mathcal{U}, V\in 
\mathcal{V} \} \, .
 \]
Next, for a stationary sequence $(Y_i)_{i \in {\mathbb{Z}}}$ of random variables with values in a Polish space $S$, define its 
 strong mixing (or $\alpha$-mixing) coefficients of order $4$  as follows:  Let  
\[
 \alpha_{\infty, 4}(n) =  \sup_{i_4 > i_3 >i_2 > i_1 \geq n }  \alpha ({\mathcal F}_{0},  \sigma (Y_{i_1}, Y_{i_2}  , Y_{i_3}, Y_{i_4}) )   \, .
 \] 
 where ${\mathcal F}_0= \sigma(Y_i, i \leq 0 )$.  As, page 146 in \cite{MPU}, these coefficients can be rewritten in the following form:  Let $B_1$ be the class of  measurable functions   from $S^4$ to ${\mathbb R} $ and bounded by one. Then
 \[
  \alpha_{\infty, 4}(n)  = \frac{1}{4}   \sup_{f \in B_1}  \sup_{i_4 > i_3 >i_2 > i_1 \geq n }   \big   \Vert  \BBE ( f( Y_{i_1}, Y_{i_2}  , Y_{i_3}, Y_{i_4})  | {\mathcal F}_0 ) -  \BBE ( f( Y_{i_1}, Y_{i_2}  , Y_{i_3}, Y_{i_4}) )  \big \Vert_1 \, .
 \]
Hence,  an application of Item (b) of Theorem \ref{thW2} provides the following result. 
 \begin{Corollary} \label{coralphamixing}
 Let $(Y_k)_{k \in {\mathbb Z}}$ be a stationary sequence of random variables  with values in a Polish space and such that  $\sum_{k \geq 1} k^2 \alpha_{\infty, 4} (k) < \infty$.
Let $f $ be a bounded measurable numerical  function and $X_k= f(Y_k) - \BBE(f(Y_k))$. Set  $S_n = \sum_{k=1}^n X_k$. 
Then  $  W_2 (P_{S_n / {\sqrt n} } ,  G_{\sigma^2}) \ll  n^{-1/2}$.
\end{Corollary}
As mentioned in the introduction, this  results applies to the class of irreducible  aperiodic  and positively recurrent Markov  $(\xi_n)$ with an atom denoted by $A$, under the condition 
$\BBE_A ( \tau_A^4) < \infty$. Here  $\tau_A$ is the first return time in $A$ and ${\BBE}_A$ stands for the expectation under ${\BBP}_x$ for $x \in A$.

\subsection{$\alpha$-dependent sequences and $\tau$-mixing sequences}\label{taualpha} 

We start by  recalling the definition of the $\alpha$-dependence coefficients as considered in  \cite{DGM10}. 
\begin{Definition}\label{defalpha}
For any random variable $Y=(Y_1, \cdots, Y_k)$ with values in
${\mathbb R}^k$ and any $\sigma$-algebra $\F,$ let
\[
\alpha(\F, Y)= \sup_{(x_1, \ldots , x_k) \in {\mathbb R}^k}
\left \| \BBE \Big(\prod_{j=1}^k (\I_{Y_j \leq x_j})^{(0)} \Big | \F \Big)^{(0)} \right\|_1.
\]
For the sequence ${\bf Y}=(Y_i)_{i \in {\mathbb Z}}$, let \begin{equation}
\label{defalpha} \alpha_{k, {\bf Y}}(0) =1 \text{ and }\alpha_{k, {\bf Y}}(n) = \max_{1 \leq l \leq
k} \ \sup_{ n\leq i_1\leq \ldots \leq i_l} \alpha(\F_0,
(Y_{i_1}, \ldots, Y_{i_l})) \text{ for $n>0$}\, ,
\end{equation}
where $\F_0 = \sigma ( Y_i, i \leq 0)$. 
\end{Definition}
Item (b) of Theorem \ref{thW2} together with  equality (A.4) in \cite{DR08} (with $\BBE ( \cdot | \F_0)$ instead of $\BBE$)  provide the following result.
\begin{Corollary} \label{coralpha}
Let $f $ be a bounded variation (BV)  function and $X_k= f(Y_k) - \BBE(f(Y_k))$ where $(Y_k)_{k \in {\mathbb Z}}$ is a stationary sequence of real-valued random variables. Let $S_n = \sum_{k=1}^n X_k$. If
$\sum_{k \geq 1} k^2 \alpha_{4, {\bf Y}} (k) < \infty$,
then  $  W_2 (P_{S_n / {\sqrt n} } ,  G_{\sigma^2}) \ll  n^{-1/2}$.  
\end{Corollary}
From this result, we can derive rates in the CLT for the partial sums associated with  BV observables of the LSV map. More precisely, for  $\gamma \in ]0,1[$, let  $T_\gamma$ defined from $[0,1]$ to $[0,1]$ by
\begin{equation*}
T_\gamma(x)=
\begin{cases}
x(1+ 2^{\gamma}x^{\gamma}) \quad \text{ if $x \in [0, 1/2[$} \\
2x-1 \quad \quad \quad \ \ \text{if $x \in [1/2, 1]$} \, .%
\end{cases}
\end{equation*}
This is the so-called LSV  \cite{LSV} map with parameter $\gamma$. Recall, that there exists a unique $T_\gamma$-invariant measure $\nu_\gamma$ on $[0, 1]$, which is absolutely
continuous with respect to the Lebesgue measure with positive density denoted by $h_\gamma$.  From Corollary \ref{coralpha} above and \cite[Prop. 1.17]{DGM10}, we derive that $  W_2 (P_{S_n / {\sqrt n} } ,  G_{\sigma^2}) \ll  n^{-1/2}$ for any $\gamma < 1/4$, where $f$ is a bounded variation function and 
$S_n = \sum_{k=1}^n  ( f(T^k_{\gamma}) -\nu_\gamma (f) )$.  
\par\medskip
We now apply  Theorem \ref{thW2}  to functions of $\tau$-dependent sequences. Before stating the result, some definitions are needed.

\begin{Definition}
Let  $\eta \in ]0, 1]$, $\ell$ be a positive integer and let  $\Lambda_\eta({\mathbb R}^\ell)$ be the set of functions $f$
from ${\mathbb R}^\ell$ to ${\mathbb R}$ such that for $x = (x_1, \dots, x_\ell)$ and $y = (y_1, \dots, y_\ell)$, 
$$
  |f(x)-f(y)|\leq \frac{1}{\ell} \sum_{i=1}^{\ell} |x_i-y_i|^{\eta} \, .
$$
Define the dependence coefficients $(\tau_{\eta, \ell,{\bf Y}} (k) )_{k \geq 1}$ of the sequence $(Y_i)_{i \in {\mathbb Z}}$ by
\[
\tau_{\eta,\ell,{\bf Y}} (k)  = \max_{1 \leq j \leq \ell }    \sup_{i_j > \ldots > i_1 \geq k} \Big \| \sup_{f \in \Lambda_\eta({\mathbb R}^j) }\Big | \BBE(f(Y_{i_1}, \dots, Y_{i_j} ) |{\mathcal F}_0)-\BBE(f(Y_{i_1}, \dots, Y_{i_j} ) )\Big | \Big \|_1 \, .
\]
\end{Definition}

Examples of $\tau_\eta$-dependent sequences are given in \cite{DP05}. 

\medskip

Let $(Y_k)_{k \in {\mathbb Z}}$ be a stationary sequence of real-valued random variables and $f$ be a bounded and $\eta$-H\"older function, with $\eta \in ]0,1]$.   Define $X_k = f(Y_k) - \BBE (f(Y_k))$. Then, for any positive integers $p,q$ and $k$, $\theta_{X, p,q } (k) \leq C \tau_{\eta, p,{\bf Y}} (k)  $ where $C$ is a positive constant depending only on $p$, $q$ and $\Vert f \Vert_{\infty}$. Hence the following result holds. 

\begin{Corollary} \label{cortaufirst}
Let $f$ be a bounded and $\eta$-H\"older function with $\eta \in ]0,1]$ and $X_k= f(Y_k) - \BBE(f(Y_k))$ where $(Y_k)_{k \in {\mathbb Z}}$ is a stationary sequence of real-valued random variables. Let $S_n = \sum_{k=1}^n X_k$. If
$\sum_{k \geq 1} k^2 \tau_{\eta, 4,{\bf Y}} (k)  < \infty$, 
then  $ W_2 (P_{S_n / {\sqrt n} } ,  G_{\sigma^2}) \ll  n^{-1/2}$.  
\end{Corollary}
From this result, we can derive rates in the CLT for the partial sums associated with H\"older functions of the LSV map above.  Starting from Corollary \ref{cortaufirst} and taking into account \cite[Prop. 5.3 and Inequality (4.2)]{DM15}, we derive that if $ \gamma <1/4$, then   $  W_2 (P_{S_n / {\sqrt n} } ,  G_{\sigma^2}) \ll  n^{-1/2}$, where $S_n = \sum_{k=1}^n  ( f(T^k_{\gamma}) -\nu_\gamma (f) ) $ and $f$ is  an $\eta$-H\"older observable with $\eta \in ]0,1]$.

\medskip

We now define another class of functions which are well adapted to  $\tau$-dependence.

\begin{Definition}
Let $c$ be any concave function from ${\mathbb R}^+$ to ${\mathbb R}^+$, with $c(0)=0$. Let ${\mathcal L}_c$
be the set of functions $g$ from ${\mathbb R}$ to ${\mathbb R}$ such that
$$
 |g(x)-g(y)| \leq K c(|x-y|), \quad \text{for some positive $K$.}
$$
\end{Definition}

Let $g \in {\mathcal L}_c$ and $X_k= g(Y_k) - \BBE(g(Y_k))$ where $(Y_k)_{k \in {\mathbb Z}}$ be a stationary sequence of bounded real-valued random variables.  
Then, for any positive integers $\ell$ and $k$,  $\tau_{1,\ell,{\bf X}} (k ) \leq K c (  \tau_{1,\ell,{\bf Y}} (k )) $. As a consequence of Corollary \ref{cortaufirst}, the following result holds:

\begin{Corollary} \label{cortau}
Let $g \in {\mathcal L}_c$ and  $X_k= g(Y_k) - \BBE(g(Y_k))$
where $(Y_k)_{k \in {\mathbb Z}}$ is a stationary sequence of bounded real-valued random variables such that $ \tau_{1,4,{\bf Y}} (k )) = O (\rho^k)$ for some $\rho$ in $]0,1[$. Let $S_n = \sum_{i=1}^n X_i$. If 
\[
\int_0^1 \frac{ (\log t)^2}{t} c(t) dt < \infty \, , 
\]
then $W_2 (P_{S_n / {\sqrt n} } ,  G_{\sigma^2}) \ll  n^{-1/2}$. 
\end{Corollary}

Corollary \ref{cortau} applies in particular to $X_k= g (T^k) - \nu(g) $ where $T$ is a map from $[0,1]$ to $[0,1]$ that can be modelled by a Young tower with exponential tails of the return times and $\nu$ is the usual invariant measure (see Section 4 in \cite{DM15} adapted to the case of exponential tails of the return times). 

\subsection{Symmetric random walk on the circle} \label{sectcircle}
$\quad \;$ Let $K$ be the Markov kernel defined by 
$Kf (x) =( f (x+a) + f (x-a) )/2$ 
on the torus $\R/\Z$, with $a$
irrational in $[0,1]$. The Lebesgue-Haar  measure $m$  is the unique
probability which is invariant by $K$. Let $(\xi_i)_{i\in \Z}$ be
the stationary Markov chain with transition kernel $K$  and
invariant distribution $m$. For $f \in {\mathbb L}^2(m)$, let
\begin{equation}\label{defSnf}
X_k=f(\xi_k)-m(f) \, .
\end{equation}
This example has been considered by Derriennic and Lin \cite{DL01} who showed  that the central
limit theorem holds with the normalization $\sqrt{n}$  as soon as \begin{equation}\label{Paroux}
 \sum_{k\in{\Z}^*}   \frac{|\hat f (k)|^2}{d (ka, {\Z})^2} < \infty
 \, ,
\end{equation}
where $\hat f(k)$ are the Fourier coefficients of $f$
and $d (ka, {\Z})=\min_{i \in \Z} |ka-i|$. The aim in this section is to give additional conditions on $f$ and on the
properties of the irrational number $a$ ensuring rates of convergence in the CLT. Let us then introduce the following definition: $a$ is said to be \textit{badly approximable in the weak sense} by rationals if for any positive $\varepsilon$, 
\beq\label{badlyweak} \text{ the inequality $d(ka, \Z) < |k|^{-1-\varepsilon}$ has only finitely many
solutions for $k \in \Z^*$.}
\eeq

\noindent From Roth's theorem the algebraic numbers are badly
approximable in the weak sense (cf. Schmidt \cite{Sc80}). Note also that the set of badly
approximable numbers in $[0,1]$ has Lebesgue measure $1$.

An application of Theorem \ref{thW2}  together with Lemma 5.2 in \cite{DR08} and their inequality (5.18) give the following corollary.

\begin{Corollary}\label{circle}
Let $X_k$ be defined by \eqref{defSnf}. Suppose that the irrational number $a$ satisfies \eqref{badlyweak}. Assume that for some positive $\varepsilon$, 
\begin{equation*} \label{condirra}
\sup_{k\not= 0} |k|^{ 6 +\varepsilon} |\hat f (k)| < \infty \, .
\end{equation*}
Then $  W_2 (P_{S_n / {\sqrt n} } ,  G_{\sigma^2}) \ll  n^{-1/2}$.  
\end{Corollary}

\section{Proofs} \label{section4}

\setcounter{equation}{0}

\subsection{Proof of Theorem \ref{thW2}} \label{subproofTh}

Assume first that $\sigma^2 =0$. In this case $G_{\sigma^2} = \delta_0$ and
\[
\BBE \big ( W_2^2 (P_{S_n / {\sqrt n} |{\mathcal F}_0 } ,  \delta_0) \big )  = n^{-1} \BBE (S_n^2) -\sigma^2 
 = - 2 n^{-1}  \sum_{k \geq 1} (k \wedge n ) {\rm Cov} (X_0,X_k) , 
\]
which, combined with \eqref{UpperBoundCov1}, shows that the upper bounds (a) and (b) hold.
\par\smallskip
We turn now to the case $\sigma^2 >0$.  Let $\delta$ be a random variable with uniform distribution over $[0,1]$ independent of $(X_k)_{k  \in {\mathbb Z}}$. Define 
${\cal G}_\ell = \sigma ((X_i)_{i \leq \ell}, \delta) $ and ${\cal G}_\infty = \sigma ((X_{i })_{i \in {\mathbb Z}}, \delta) $. Define also  the conditional expectation operator $\BBE_0$ by  $\BBE_0 ( \cdot )= \BBE ( \cdot | {{\mathcal  G}_0 } ) $.

\smallskip

In what follows $(Y_k)_{k \geq 1}$ will be a sequence of iid random variables  independent of ${\cal G}_\infty$.  In case of Item (a), their common law will be the normal law  $ {\cal N}(0, \sigma^2)$ whereas in case of Item (b), we will have to prescribe also their third moment as it is described below. 

Let  $\beta_3$ be a fixed real number. Let $Z$ be a r.v.  with distribution $ {\cal N}(0, \sigma^2/2)$. There exists a random variable 
$B$ independent of $Z$, taking only 2 values and such that  $Y= Z+B$ satisfies 
\beq \label{momentY}
{\mathbb E}(Y)=0 \, , \, {\mathbb E}(Y^2)= \sigma^2 \text{ and }
{\mathbb E}(Y^3)=\beta_3 \, .
\eeq
 We refer to Lemma 5.1 in \cite{DR08} for more details. 
For the proof of Item (b), 
\beq \label{defbeta3}
\beta_3= \BBE (X_0^3) + 3 \sum_{i \geq 1 }  \{ \BBE ( X_0^2 X_i)  + \BBE ( X_0 X^2_i)  \} + 6 \sum_{u \geq 1} \sum_{v \geq u+1} \BBE ( X_0 X_uX_{v}) \, ,
\eeq
which is the limit of $n^{-1}\BBE (S_n^3) $, as $n \rightarrow \infty$, under the conditions of Item (b) of Theorem \ref{thW2}.

Let $(Z_k)_{k \geq 1}$ be a sequence of independent r.v.'s distributed as $Z$ and let $(B_k)_{k \geq 1}$ be a sequence of independent r.v.'s distributed as $B$ and independent of $(Z_k)_{k \geq 1}$. Suppose furthermore that the sequence $(Z_k,B_k)_{k \geq 1}$  is independent of ${\cal G}_\infty$.  For any $k \geq 1$, set $Y_k = Z_k+B_k$.  

Next, in case of both items, we define $T_n = Y_1 + Y_2 + \cdots + Y_n$.  Note that 
\[
W_2 ( P_{S_n / {\sqrt n}|{\mathcal F}_0 } , G_{\sigma^2} )  \leq  W_2 ( P_{S_n / {\sqrt n}|{\mathcal F}_0 } , P_{T_n / {\sqrt n}}  )  + W_2 ( P_{T_n / {\sqrt n}} , G_{\sigma^2}  ) \, .  
\] 
According to Theorem 4.1 in \cite{Rio09}, since $Y \in {\mathbb L}^4$, $W_2 ( P_{T_n / {\sqrt n}} , G_{\sigma^2}  ) \ll n^{-1/2}$. Since  $P_{S_n / {\sqrt n}|{\mathcal F}_0} = P_{S_n / {\sqrt n}|{\mathcal G}_0}$, the theorem will follow if one can prove that the upper bounds (a) and (b) still hold with $P_{S_n / {\sqrt n}|{\mathcal G}_0 }$ replacing $P_{S_n / {\sqrt n}|{\mathcal F}_0 }$. 
With this aim, we shall apply Lemma 5.1 in \cite{MR12}.  We start by introducing some notations.  Let $\Lambda_2$ be the class of
real functions $f$ which are continuously differentiable
and such that $|f'(x) - f'(y) |\leq | x - y |$  for any  $(x,y) \in {\mathbb R} \times {\mathbb R}$. Let also $\Lambda_2 (E)$ be  the set of measurable functions $f:{\mathbb R}
\times  E \rightarrow {\mathbb R}$ wrt the $\sigma$-fields 
${\mathcal L} ( {\mathbb R} \times E) $ and ${\mathcal B} ({\mathbb R})$, 
such that $f( \cdot, w) \in \Lambda_2$ and $f(0,w)=f'(0,w)=0$  for any 
$w \in E$. Next, let $W= ( (X_i)_{ i \in {\Z}^-}, \delta)$ and $E = {\R}^{{\Z}^-} \times [0,1]$. According to \cite[Lemma 5.1]{MR12}, and denoting by $N$   a ${\mathcal N}(0,\sigma^2)$-distributed random variable, independent of  all the above  sequences (so independent of $(X_k,Y_k)_k$), the upper bound (a) will follow  if one can prove that 
\beq \label{W2condp2-Itema}
\sup_{f \in \Lambda_2 ( E) }  \BBE \big (  f ( S_n + N  , W   ) -    f ( T_n + N  , W   )  \big )  \ll  \sqrt{n} \Big (  1 + \sum_{k\geq 1}  ( k \wedge \sqrt{n} )  \theta_{X,2,2} (k) \Big )  \, ,
\eeq
whereas the upper bound (b) will follow  if
\beq \label{W2condp2}
\sup_{f \in \Lambda_2 ( E) }  \BBE \big (  f ( S_n + N  , W   ) -    f ( T_n + N  , W   )  \big )  \ll  1 + \sum_{k\geq 1} k ( k \wedge \sqrt{n} )  \theta_{X,4,4} (k) \, .
\eeq
In what follows to soothe the notations we omit the subscripts for the coefficients $\theta (k)$. 

\smallskip

\noindent {\it Proof of Item (a).}  We shall apply the  Lindeberg method. 
Let us first introduce some notations.
\begin{Notation}\label{not21-a}
Set 
$f_{n-k} (x ) 
% := f_{n-k} (x, W ) &
 =  {\mathbb E}_0 ( f (x + N + T_n - T_k , W)  )$.
% \\
%&  = \int_{\mathbb R}  f (x + t , W)  P_{N + T_n - T_k} (dt) \, .
%\end{align*}
\end{Notation}
Notice that 
\[
f_{n-k} (x )  =  \int_{\mathbb R}  f (x - t , W) \varphi_{\sigma^2 (n-k+1)} (  t)  dt  \, , 
\]
where $\varphi_{t^2}$ is the density of a ${\mathcal N}(0,t^2)$.  Hence, according to Lemma 6.1 in \cite{DMR09}, 
\beq \label{boundderiv}
\Vert  f_{n-k}^{(i) }\Vert_{\infty} := b_i  \ll (n-k+1)^{(2-i)/2} \, .
\eeq
Since the sequence $(N,(Y_i)_{i \geq 1})$ is independent of the sequence $(\XZ,W)$,
\begin{equation}\label{sumdeltaalter}
{\mathbb E} ( f(S_n +  N, W) -f (T_n +  N,W) ) = \sum_{k=1}^n {\mathbb E} ( f_{n-k} (S_{k-1} +X_k)  - f_{n-k} (S_{k-1} +Y_k)  ) \, .
\end{equation}
By the Taylor formula at order 3 and using \eqref{boundderiv},  we get 
\begin{equation}\label{sumdelta-1}
\Big | \BBE \big ( f_{n-k} (S_{k-1} + Y_k) -f_{n-k} ( S_{k-1}) - \frac{\sigma^2}{2} f''_{n-k} ( S_{k-1}) \big ) \Big | \leq
C (n-k+1)^{-1/2} \, .
\end{equation}
Similarly 
\begin{equation}\label{sumdelta-2}
 \big | \BBE (  f_{n-k}( S_k )   -f_{n-k} (S_{k-1} )    - f_{n-k}^\prime (S_{k-1}) X_k -
{1\over 2}f_{n-k}'' (S_{k-1}) X_k^2  ) \big | \leq
C (n-k+1)^{-1/2} \, .
\end{equation}
Now we control the second order term.
Let \beq \label{defgamma} \Gamma_{k}(k,i) = f_{n-k}'' (S_{k-i}) -f_{n-k}''  (S_{k-i-1}) \, .\eeq
Clearly
$$
f_{n-k}'' (S_{k-1})  X_k^2 = \sum_{i=1}^{[\sqrt{k}]-1}
\Gamma_{k}(k,i) X_k^2  +f_{n-k}'' (S_{k-[\sqrt k]})X_k^2 \, .
$$
Since
$|\Gamma_{k}(k,i) | \leq b_3 |X_{k-i}|$, by stationarity we get that for any $i \leq k-1$,
$$
\big |\cov(\Gamma_{k}(k,i) , X^2_k   ) \big | \leq  b_3 \Vert X_0 \big (\BBE_0 (X^2_i  )- \BBE (X^2_i ) \big ) \Vert_1  \ll  (n-k+1)^{-1/2} \theta(i)\, .
$$
Since $\Vert f_{n-k}'' \Vert_{\infty} \leq b_2 $ a.s., we also get
by stationarity that
$$
\big |\cov ( f_{n-k}'' (S_{k-[\sqrt k]}), X^2_k) \big | \leq  b_2 \Vert \BBE_0 (X_{[\sqrt k]}^2)-  \BBE (X_{[\sqrt k]}^2)\Vert_1 \ll  \theta([\sqrt k]) \, .
$$
Starting from \eqref{sumdelta-2}, it follows that 
\begin{multline} \label{dt5}
\big | \BBE (  f_{n-k}( S_k )   -f_{n-k} (S_{k-1} )    - f_{n-k}^\prime (S_{k-1}) X_k) -
{1\over 2} \BBE (f_{n-k}'' (S_{k-1})) \BBE (X_k^2) \big | \\ \ll    \theta([\sqrt k]) + (n-k+1)^{-1/2}  \big (1+  \sum_{i=1}^{[\sqrt{k}]}\theta(i) \big ) 
  \,  .
\end{multline}
Starting from \eqref{sumdeltaalter} and taking into account \eqref{sumdelta-1} and  \eqref{dt5}  we derive that 
\begin{multline} \label{sumdelta-3}
 \big | {\mathbb E} ( f(S_n +  Y, W) -f (T_n +  Y,W) )  \big | \\
   \ll \sqrt{n}  \Big ( 1 + \sum_{i=1}^{[\sqrt n]} \theta (i)  \Big )  + \Big |    \sum_{k=1}^n  \Big \{  \BBE ( f'_{n-k} (S_{k-1} )  X_k ) -  \BBE (f_{n-k}'' (S_{k-1})) \sum_{j \geq 1} \BBE(X_0 X_j)  \Big \}  \Big | 
  \,  .
\end{multline}
To give now an estimate of $\BBE (f'_{n-k} ( S_{k-1})X_k)$, we write
\[
f'_{n-k} (S_{k-1} ) = f'_{n-k} (0) + \sum_{i=1}^{k-1} ( f'_{n-k}(S_{k-i}) - f'_{n-k}( S_{k-i-1}) ) \, .
\]
Hence
\beq\label{dt6}
\BBE (f'_{n-k} ( S_{k-1})X_k)  =   \sum_{i=1}^{k-1} {\rm Cov} \big (f'_{n-k}(S_{k-i}) - f'_{n-k}( S_{k-i-1} ) , X_k \big ) +
\BBE ( f'_{n-k} (0) X_k ) \, .
\eeq
Now $f'_{n-k} (0) $ is a ${\cal G}_0$-measurable random variable. Since $f \in \Lambda_2(E) $ then $f'(0,w)=0$ and $f' (\cdot, w) $ is $1$-Lipschitz. Therefore  
$$
|f'_{n-k} (0)| \leq  \int_{\mathbb R} | f'(u,W) -f'(0,W) |   \varphi_{\sigma^2 (n-k+1 ) } (-u) du \leq  \sigma \sqrt{ n-k+1 } \, \text{ a.s.}
$$
It follows that 
\beq \label{dt7}
\sum_{k=1}^n  \big | {\mathbb E} ( f'_{n-k} (0) X_k )  \big |  \ll  \sum_{k=1}^n \sqrt{ n-k+1} \Vert \BBE_0 (X_k) \Vert_1 \ll \sqrt{n} \sum_{k=1}^n  \theta(k) \, .
\eeq
We give now  an estimate of $\sum_{i=1}^{k-1} {\rm Cov} \big (f'_{n-k}(S_{k-i}) - f'_{n-k}( S_{k-i-1} ) , X_k \big )$. Using the 
stationarity and noting that $|f_{n-k}^\prime (S_{k-i} ) -f_{n-k}^\prime (S_{k-i-1})| \leq b_2|X_{k-i}|$,
we have
\[
|\cov (f_{n-k}^\prime (S_{k-i} ) -f_{n-k}^\prime (S_{k-i-1} ),  X_k)| \leq b_2 M  \Vert  \BBE_0 ( X_i) \Vert_1 \ll  \theta(i) \, .
\]
Hence 
\begin{align} \label{dt8}
 \sum_{k=1}^n \sum_{i=[\sqrt k]}^k |\cov (f_{n-k}^\prime (S_{k-i} ) & -f_{n-k}^\prime (S_{k-i-1} ),  X_k)| \nonumber  \\ &  \ll  \sum_{i =1}^n (i \wedge  \sqrt{n})^2\theta(i) 
  \ll \sqrt{n}  \sum_{i \geq 1} (i \wedge \sqrt{n})\theta(i)  \, .
\end{align}
From now on, we assume that $i< [{\sqrt k}]$. We first write 
\[f_{n-k}^\prime (S_{k-i} ) -f_{n-k}^\prime (S_{k-i-1} ) =f_{n-k}'' (S_{k-i-1})  X_{k-i} +
  R_{k,i} \, ,
\]
where $R_{k,i}$ is ${\cal F}_{k-i}$-measurable and
$|R_{k,i}| \leq b_3 X_{k-i}^2 /2$. Hence, by stationarity,
\[
|\cov (R_{k,i},  X_k) | \leq b_3 \Vert X^2_0 \BBE_0 (  X_i) \Vert_1 /2 \ll (n-k+1)^{-1/2} \theta(i) \, .
\]
implying that 
\begin{equation} \label{dt10}
\sum_{k=1}^n \sum_{i=1}^{[\sqrt k]}|\cov (R_{k,i},  X_k) | \ll  \sqrt{n} \sum_{i = 1}^{[\sqrt n]}  \theta(i) \, .
\end{equation}
In order to estimate the term $\BBE (f''_{n-k} (S_{k-i-1})X_{k-i} X_k)$,
we introduce the decomposition below:
$$
f''_{n-k} (S_{k-i-1}) = \sum_{\ell=1}^{i-1} ( f''_{n-k}
(S_{k-i-\ell}) - f''_{n-k} (S_{k-i-\ell-1}) )  + f''_{n-k} (S_{k-2i}) \,  ,
$$
where by convention we set $S_p = 0$ if $p \leq 0$. 
For any $\ell \in \{ 1, \cdots, i-1 \}$, by using the notation (\ref{defgamma}) and the stationarity, we get that
\[
|\cov (\Gamma_{k}(k,\ell+i) X_{k-i} , X_k )| \leq b_3 \Vert X_{- \ell} X_0 \BBE_0 ( X_i ) \Vert_1  \ll (n-k+1)^{-1/2} \theta(i)  \, .
\]
Hence 
\begin{equation} \label{dt11}
\sum_{k=1}^n \sum_{i=1}^{[\sqrt k]}  \sum_{\ell =1}^{i-1} \cov (\Gamma_{k}(k,\ell+i) X_{k-i} , X_k )|  \ll  \sqrt n\sum_{i=1}^{[\sqrt n]}   i \theta(i) \, .
\end{equation}
As a second step, we bound up
$|\cov (f_{n-k}'' (S_{k-2i}) ,  X_{k-i}   X_k )|$.  Clearly, 
$$
f_{n-k}'' (S_{k-2i}) = \sum_{\ell =i}^{k-i -1} \Gamma_{k}(k, \ell+i) + f_{n-k}'' (0) \, .
$$
Now for any $\ell \in \{ i, \cdots,  (k-i-1) \}$, by  stationarity,
$$
|\cov (\Gamma_{k}(k, \ell+i)   , X_{k-i} X_k )| \leq  b_3 \Vert X_{- \ell} \big (\BBE_{- \ell} (X_0 X_i )
- \BBE ( X_0  X_i ) \big ) \Vert_1 \ll (n-k+1)^{-1/2} \theta(\ell)   \, .
$$
Hence 
\begin{equation} \label{dt12}
\sum_{k=1}^n \sum_{i=1}^{[\sqrt k]}  \sum_{\ell=i}^{k-i -1}   |\cov (\Gamma_{k}( k, \ell +i)  , X_{k-i} X_k )|  \ll  \sqrt n \sum_{\ell=1}^n (\ell \wedge \sqrt{n})\theta( \ell)   \, .
\end{equation}
Next, note that 
\[
|\cov (f_{n-k}'' (0 ), X_{k-i} X_k )| \ll b_2 \min ( \theta (k-i), \theta(i) )  \ll \theta([k/2])  \, ,
\]
implying that 
\beq \label{dt13second}
\sum_{k=1}^n \sum_{i=1}^{[\sqrt k]}  |\cov (f''_{n-k}(0 ), X_{k-i} X_k )| \ll \sum_{k=1}^n \sum_{i=1}^{[\sqrt k]}   \theta([k/2]) \ll \sqrt n \sum_{k=1}^n\theta( k)   \, .
\eeq
Taking into account the inequalities (\ref{dt7})-(\ref{dt13second}), and using that $\sum_{k \geq 1}\theta( k) < \infty$,  we get 
\begin{equation}\label{dt14} 
 \sum_{k=1}^n \Big | \BBE (f^\prime_{n-k} (  S_{k-1})  X_k )  -  \sum_{i=1}^{ [\sqrt k] }
\BBE (f_{n-k}'' (S_{k-2i})) \BBE  (X_{k-i} X_k )  \Big | \ll \sqrt n \Big ( 1 +  \sum_{\ell \geq 1} (\ell \wedge \sqrt{n})\theta( \ell)  \Big )  \, .
 \end{equation}
We handle now the quantity 
\[
A_k:= \sum_{i=1}^{ [\sqrt k] } \BBE (f''_{n-k}(S_{k-2i})) \BBE ( X_{k-i}
 X_k) - \sum_{i=1}^{\infty} \BBE (f''_{n-k}
(S_{k-1}) )  \BBE (X_{k-i}X_k) \, .
\]
We first note that by stationarity,
\[
\sum_{i\geq [\sqrt k]  +1}| \BBE (f''_{n-k}(S_{k-1}) )  \BBE (X_{k-i}X_k) |
\leq  b_2 \sum_{i\geq [\sqrt k]  +1} | \BBE( X_0 \BBE_0(X_i)) | \ll\sum_{i\geq [\sqrt k]  +1}  \theta(i)  \, .
\]
Hence 
\begin{equation}\label{dt14-2} 
\sum_{k=1}^n \sum_{i\geq [\sqrt k]  +1}| \BBE (f''_{n-k}(S_{k-1}) )  \BBE (X_{k-i}X_k) | \ll \sum_{i\geq 1 } (i \wedge \sqrt n)^2   \theta(i)  \ll \sqrt n  \sum_{i\geq 1 } (i \wedge \sqrt n)   \theta(i)  \, .
 \end{equation}
On another hand, we write
\[
\BBE (f''_{n-k}(S_{k-1}) -f_{n-k}'' (S_{k-2i})) \BBE ( X_{k-i}
 X_k ) = \sum_{\ell =1}^{2i-1} \BBE (\Gamma_{k}(k,\ell )  \BBE (
X_0 \BBE_0( X_i) )\, .
\]
Therefore
\[
\sum_{i=1}^{[\sqrt k]}| \BBE (f''_{n-k}(S_{k-1}) -f_{n-k}'' (S_{k-2i})) \BBE ( X_{k-i} X_k ) | 
 \leq   (n-k+1)^{-1/2}  \sum_{i=1}^{[\sqrt k]}  i  \theta(i) \, , 
\]
implying that 
\begin{equation} \label{dt18}
\sum_{k=1}^n \sum_{i=1}^{[\sqrt k]}| \BBE (f''_{n-k}(S_{k-1}) -f_{n-k}'' (S_{k-2i})) \BBE ( X_{k-i}  X_k ) |  \ll  {\sqrt n}
\sum_{i=1}^{[\sqrt n]}  i  \theta(i)   \, .
\end{equation}
Hence (\ref{dt14-2}) and (\ref{dt18}) entail that
\begin{equation} \label{dt19}
\sum_{k=1}^n |A_k|  \ll \sqrt n  \sum_{i\geq 1 } (i \wedge \sqrt n)   \theta(i)  \, .  
\end{equation}
The estimates \eqref{dt14}  and  \eqref{dt19} yield to  
\begin{equation}\label{dt20} 
\sum_{k=1}^n  \Big  |\BBE (f^\prime_{n-k} (  S_{k-1})  X_k )  -  \sum_{i=1}^{\infty} \BBE (f''_{n-k}
(S_{k-1}) )  \BBE (X_{0}X_i)  \Big  | \ll \sqrt n \Big ( 1 +  \sum_{\ell \geq 1} (\ell \wedge {\sqrt n})\theta( \ell)  \Big )  \, .
 \end{equation}
Taking into account the estimates  \eqref{sumdelta-3} and  \eqref{dt20}, Item  (a) follows.

\medskip

\noindent {\it Proof of Item (b).} Recall that in this case the iid random variables $(Y_k)_{k \geq 1}$  have their first three moments defined by \eqref{momentY}  and \eqref{defbeta3}

\begin{Notation} \label{TildeX} For any integer $k \geq 0$, let 
$\X_k = X_k - \BBE_0 (X_k)$ and $\S_k = S_k - \BBE_0 (S_k)$, with the convention $S_0=0$.
\end{Notation}

Note that, since we assume that $ \sum_{j\geq 1} j \theta (j) < \infty$, 
\[
\Vert \BBE_0 (S_n) \Vert_2^2 \leq 2 \sum_{i=1}^n \sum_{j =i }^n  \big | \BBE ( \BBE_0 (X_i) \BBE_0 (X_j)  )  \big |  \leq 2 M \sum_{j=1}^n j \theta (j) \ll 1  \, .
\]
Therefore, using that $f'(0,W) =0$ and that $ | f'(x,W) -  f'(y,W)  | \leq |x-y|$, we infer that  to prove \eqref{W2condp2}, it is enough to show  that for any $f \in \Lambda_2(E) $ and any positive $n$, 
\beq \label{aim1thmain}
\sup_{f \in \Lambda_2 ( E) }  \BBE \big (  f ( \S_n + N  , W   ) -    f ( T_n + NN  , W   )  \big )  \ll  1 + \sum_{k\geq 1} k ( k \wedge \sqrt{n} )  \theta (k) \, .
\eeq
This will be done by using again the Lindeberg method. Let us introduce some additional notations.
\begin{Notation}\label{not21}
%Set 
%$f_{n-k} (x ) 
%% := f_{n-k} (x, W ) &
% =  {\mathbb E}_0 ( f (x + Y + T_n - T_k , W)  )$.
%% \\
%%&  = \int_{\mathbb R}  f (x + t , W)  P_{Y + T_n - T_k} (dt) \, .
%%\end{align*}
For any positive integer $k$, 
let $\Delta_{n,k} =  f_{n-k} ({\tilde S}_{k-1} + {\tilde X}_k) - f_{n-k} ( {\tilde S}_{k-1} + Y_k)$ where $f_{n-k}$ is defined in Notation \ref{not21-a}.
\end{Notation}
All along the proof, the following lemma will be used (the proof is postponed to the Appendix and is based on the fact that the common distribution of the random variables $(Y_k)_{k \geq 1}$ is smooth).  
\begin{Lemma} \label{lmacrucial} Let $f \in \Lambda_2(E)$. 
\begin{enumerate}
\item[1.] For any $i \geq 2$, there exists a positive constant  $\kappa_1$ depending on $\sigma^2$ and  $i$ and such that $\Vert  f_{n-k}^{(i) }\Vert_{\infty} \leq \kappa_1(n-k+1)^{(2-i)/2}$.

\item[2.] Assume that $\sum_{k \geq 1} k \theta_{X,3,4} (k) < \infty$. Then, for any $i \geq 2$, 
there exists a constant $\kappa_2>0$ depending on $\sigma^2$ and  $i$ 
such that, for any  integer $\ell >0$, 
\[
\big | \BBE (  f_{n-k}^{(i) } ({ \tilde S}_{\ell -1}  )  )  \big | \leq \kappa_2(n-k+1)^{(1-i)/2}  + 
\kappa_2(n-k+ \ell)^{(2-i)/2} \, .
\]
\end{enumerate}
\end{Lemma}
\begin{Remark}If $(X_k)_{k \in {\mathbb Z}}$ is a stationary sequence of  martingale differences, Item 2. is valid under the condition $\sum_{k \geq 1}  \theta_{X,2,3} (k) < \infty$. 
\end{Remark}
\par
Since the sequence $(N,(Y_i)_{i \geq 1})$  is independent of  $( \XN , W)$,
\begin{equation}\label{sumdelta}
{\mathbb E} ( f ({ \tilde S}_n +  N,W) - f (T_n +  N,W) ) = \sum_{k=1}^n {\mathbb E} ( \Delta_{n,k} )  \, .
\end{equation}
Next the functions $f_{n-k}$ are $C^\infty$. Consequently, from the Taylor integral formula at order $5$,
\beq \label{firstdecomp}
 \Delta_{n,k} =  \sum_{j=1}^4 \frac{1}{j!} \, f^{(j)}_{n-k} ({ \tilde S}_{k-1}) ({ \tilde X}_k^j - Y_k^j)  +R_{n,k} \, ,
\eeq
with
\[
R_{n,k} =  \frac{1}{24}   { \tilde X}_k^5  \int_0^1 (1-s)^4 f_{n-k}^{(5)} ({ \tilde S}_{k-1}  + s { \tilde X}_k) ds \\
 -  \frac{1}{24}   Y_k^5  \int_0^1 (1-s)^4 f_{n-k}^{(5)} ({ \tilde S}_{k-1}  + s Y_k) ds  \, .
\]
Taking into account the fact that $\Vert X_{k} \Vert_{\infty}  \leq M $ and Item 1 of Lemma \ref{lmacrucial}, we derive that 
\[
\Vert R_{n,k} \Vert_1   \ll  ( M^5  + \BBE ( |Y_1|^5 )  )  \Vert  f_{n-k}^{(5) }\Vert_{\infty}    
 \ll    (n-k+1)^{-3/2}  \, .
\]
Therefore, 
\begin{equation}\label{reste1}  
\sum_{k\in [1,n]} \Vert R_{n,k} \Vert_1  \ll 1 \, .
\end{equation}
Let $\beta_4= \sigma^2= \BBE(Y_k^2)$ and $\beta_4= \BBE(Y_k^4)$. Since the sequence $(Y_i)_{i \geq 1}$ is independent of the sequence $(X_i)_{i \geq 1}$,  
\begin{align} \label{linddec1}
\BBE  (   \Delta_{n,k}   -  R_{n,k}  ) &  = \BBE \Bigl  ( f'_{n-k} ({\tilde S}_{k-1}) \X_k    +  \sum_{\ell =2}^4  \frac{1}{\ell ! } f^{(\ell)}_{n-k}  ({\tilde S}_{k-1}) ( \X_k^\ell - \beta_\ell) \Bigr)  \nonumber  \\
& =  \BBE \Bigl( f'_{n-k} ({\tilde S}_{k-1}) \X_k    +     \sum_{\ell =2}^4  \frac{1}{\ell ! } f^{(\ell)}_{n-k}  ({\tilde S}_{k-1}) ( X_k^\ell - \beta_\ell)  \Bigr)   +   {\tilde B}_{n,k} 
\nonumber  \\
& :=  \BBE \Bigl( \Delta_{n,k}^{(1)} + { \frac{1}{2} }  \Delta_{n,k}^{(2)} +  { \frac{1}{6} }    \Delta_{n,k}^{(3)}  + \frac{1}{24}   \Delta_{n,k}^{(4)}  \Bigr) 
 +  {\tilde B}_{n,k}\, .
\end{align}
Using  Item 1 of Lemma \ref{lmacrucial},  we first notice that
\beq \label{tildeB}
\sum_{k\in [1,n]} \vert  {\tilde B}_{n,k} \vert  \ll  \sum_{k\in [1,n]} \Vert \BBE_0 (X_k) \Vert_1 \ll 1 .
\eeq
Next we develop the first four terms in the right-hand side of the  decomposition \eqref{linddec1} with the help of the Lindeberg method.  From now on, to soothe the notation, we shall omit most of the time the index $n$ in all the $\Delta_{n,k}^{(i)} $  and the related quantities,  and then rather write  $\Delta_{k}^{(i)} $. Let us start with the term 
$  \Delta_{k}^{(4)} $. Using  Item 2 of Lemma \ref{lmacrucial},  note first  that
\begin{equation} \label{T4correction1}
 \sum_{k=1}^n \big | \BBE  (     f_{n-k}^{(4) }    ({\tilde S}_{k-1})  ) ( \BBE (X_k^4)  + \beta_4  ) \big |    \ll  \sum_{k=1}^n  \Big (  \frac{1}{(n-k+1)^{3/2}}+  \frac{1}{n} \Big ) \ll 1 \, .
 \end{equation}
 Next, we write 
\begin{multline*} 
    f_{n-k}^{(4) }    ({\tilde S}_{k-1})   ( X_k^4 -  \BBE (X_k^4) ) =   f_{n-k}^{(4) }    (0)   ( X_k^4 -  \BBE (X_k^4) )  \\ +  \sum_{i=1}^{k-1} \big  (   f_{n-k}^{(4) }  ({\tilde S}_{k-i})  -   f_{n-k}^{(4) }  ({\tilde S}_{k-i -1 })   \big )  (X_k^4 - \BBE (X_k^4) )  \, .
\end{multline*} 
By Item 1 of Lemma \ref{lmacrucial} we get 
\begin{equation} \label{correction3}
 \sum_{k=1}^n   \big |   {\rm Cov}   (     f_{n-k}^{(4) }    (0)  , X_k^4  ) \big |     \ll 
  \sum_{k=1}^n    (n-k+1)^{-1}  \theta(k)  \ll    \sum_{k=1}^n    \theta(k)   \, ,
\end{equation}
and 
\begin{multline} \label{correction1}
 \sum_{k=1}^n \sum_{i=1}^{k-1}  \big |  {\rm Cov}   \big (   f_{n-k}^{(4) }  ({\tilde S}_{k-i})  -   f_{n-k}^{(4) }  ({\tilde S}_{k-i -1 })  , X_k^4    \big ) \big |   \\  \ll
M^5  \sum_{k=1}^n  (n-k+1)^{-3/2} \sum_{i=1}^k \theta(i)  \ll  \sum_{i=1}^n \theta(i)   \, .
\end{multline}
Taking into account \eqref{T4correction1}, \eqref{correction3}, \eqref{correction1} and the fact that $\sum_{k \geq 1}\theta (k) < \infty$, it follows that 
\begin{equation} \label{newtermorder4}
 \sum_{k=1}^n   \big | \BBE (\Delta_{n, k}^{(4)} )  \big | \ll 1 \, .
 \end{equation}
Now,  concerning the first term in the right-hand side of \eqref{linddec1}, letting
${\ell_k} = [k/2]$,   we write
\begin{align} \label{linddec2tilde}
 \BBE ( \Delta_{k}^{(1)}) & =  \BBE ( f'_{n-k} (\S_{k- {\ell_k} -1}) \X_k  )  + \sum_{i=1}^{\ell_k}  \BBE \big ( \{  f'_{n-k} (\S_{k-i}) - f'_{n-k} (\S_{k-i- 1})  \} \X_k   \big )   \nonumber \\
&  =    \BBE  ({\tilde \Delta}_{k,2}^{(1)}  )   +  \frac{1}{2}   \BBE  ({\tilde \Delta}_{k,3}^{(1)} )  + \frac{1}{6}   \BBE  ({\tilde \Delta}_{k,4}^{(1)} )   + {\tilde B}_{n,k}^{(1)}     \, ,
\end{align}
where, for $j=2,3,4$, 
\[
{\tilde \Delta}_{k,j}^{(1)} = \sum_{i=1}^{\ell_k}      f^{(j)}_{n-k} (\S_{k-i- 1})    \X^{j-1}_{k-i} \X_k   \, ,      
\]
and
\[
{\tilde B}_{n,k}^{(1)}   =   \BBE ( f'_{n-k} (S_{k- {\ell_k} -1}) \X_k  )     +   \frac{1}{6}  \sum_{i=1}^{\ell_k}  \int_0^1 (1-s)^3 \BBE \big (   f_{n-k}^{(5)} (\S_{k-i-1}  + s \X_{k-i})  \X_{k-i}^4 \X_k  \big ) ds . 
\]
We start by noticing that, by Item 1 of Lemma \ref{lmacrucial}, for any  $m \geq 2 $ and any  $s$ in $[0,1]$,  
\begin{multline} \label{compwithtilde}
 \sum_{k=1}^n \sum_{i=1}^{\ell_k}   \Vert     f^{(m)}_{n-k} (\S_{k-i- 1} + s  \X_{k-i})      (   \X^{m-1} _{k-i}\X_k - X^{m-1} _{k-i} X_k )  \Vert_1 \\  \ll
M^{m-1}  \sum_{k=1}^n  (n-k+1)^{(2-m)/2}\sum_{i=1}^{\ell_k}   \big (  \Vert \BBE_0 (X_k ) \Vert_1  + \Vert  \BBE_0 (X_{k-i} )  \Vert_1 \big )  \ll  \sum_{k \geq 1} k \theta(k)  \, .
\end{multline}
On another hand, since $ f'_{n-k} (0) $ is ${\mathcal F}_0$-measurable, $\BBE (  f'_{n-k} (0) \X_k ) =0$. Therefore
\begin{multline*}
\vert \BBE ( f'_{n-k} (\S_{k- {\ell_k} -1}) \X_k  )  \vert =  \vert \BBE  \big ( \{  f'_{n-k} (\S_{k- {\ell_k} -1})  -   f'_{n-k} (0)  \} \X_k   \big )  \vert  \\
\leq   \int_0^1   \vert \BBE  \big (f''_{n-k} ( t \S_{k- {\ell_k} -1}) \S_{k- {\ell_k} -1} \X_k   \big )  \vert  dt  \leq 4  M   \Vert   f_{n-k}^{(2)}  \Vert_{\infty}   (k-{\ell_k})\theta({\ell_k})  \, .
\end{multline*}
Since $\Vert   f_{n-k}^{(2)}  \Vert_{\infty} \ll 1$ and  $\sum_{k \geq 1} k \theta(k) < \infty$, 
\beq \label{B1b1}
 \sum_{k \in [1,n]} \vert \BBE ( f'_{n-k} (\S_{k- {\ell_k} -1}) \X_k  )  \vert  \ll 1 \, .
\eeq 
Next, Item 1 of Lemma \ref{lmacrucial} implies that 
\[
 \vert    \BBE \big (   f_{n-k}^{(5)} (\S_{k-i-1}  + s \X_{k-i})  X_{k-i}^4 X_k  \big )   \vert   \leq M^4  \Vert   f_{n-k}^{(5)}  \Vert_{\infty}  \theta (i)   \ll (n-k+1)^{-3/2}  \theta(i)  \, .
\]
Hence 
\beq  \label{B1b3}
 \sum_{k =1}^n \sum_{i=1}^{\ell_k}  \vert    \BBE \big (   f_{n-k}^{(5)} (S_{k-i-1}  + s X_{k-i})  X_{k-i}^4 X_k  \big )   \vert  \ll \sum_{i=1}^{n}  \theta(i)  \ll 1 \, .
\eeq
The upper bounds  \eqref{compwithtilde}, \eqref{B1b1} and \eqref{B1b3} imply that  
\beq  \label{B1tilde}
 \sum_{k\in [1,n]} \vert {\tilde B}_{n,k}^{(1)}   \vert   \ll 1 \, .
\eeq
Next, taking into account  Item 2 of Lemma \ref{lmacrucial} and the fact that $|\BBE( X^3_{k-i}X_k ) | \leq M^3\theta (i) $, we derive that 
\[
 \vert   \BBE \{  f^{(4)}_{n-k} (\S_{k-i- 1})  \}    \BBE( X^3_{k-i}X_k )   \vert  \ll  ( (n-k+1)^{-3/2} + (n-i)^{-1}) \theta(i)  \, .
\]
Therefore
\beq \label{B1b2}
  \sum_{k=1}^n  \sum_{i=1}^{\ell_k}  \vert   \BBE \{  f^{(4)}_{n-k} (\S_{k-i- 1})  \}    \BBE( X^3_{k-i}X_k )   \vert  \ll   \sum_{i \geq 1}\theta(i)  \ll 1 \, .
\eeq
So, overall, starting from \eqref{linddec2tilde} and taking into account  \eqref{compwithtilde},  \eqref{B1tilde} and  \eqref{B1b2}  we get
\beq\label{linddec2}
 \BBE ( \Delta_{k}^{(1)}) 
  =    \BBE  (\Delta_{k,2}^{(1)} )  +  \frac{1}{2}   \BBE  (\Delta_{k,3}^{(1)} )  + \frac{1}{6}   \BBE  (\Delta_{k,4}^{(1)} )   +  A_{k,2}^{(1)}  +   \frac{1}{2}   A_{k,3}^{(1)} + B_{n,k}^{(1)}     \, ,
\eeq
where $B_{n,k}^{(1)}$ is such that   
\beq \label{BoundB1}
\sum_{k\in [1,n]} | B_{n,k}^{(1)}  |  \ll 1 \, ,
\eeq
 and the following notations have been used:  for $j= 2,3,4$, 
\begin{equation} \label{defdeltank4}
 \Delta_{k,j}^{(1)} = \sum_{i=1}^{\ell_k}     \{  f^{(j)}_{n-k} (\S_{k-i- 1})  \}   ( X_{k-i}^{j-1}X_k )^{(0)},   
A_{k,j}^{(1)} = \sum_{i=1}^{\ell_k}    \BBE \{  f^{(j)}_{n-k} (\S_{k-i- 1})  \}    \BBE( X_{k-i}^{j-1}X_k ) .
 \end{equation}
Introduce now the following additional notations.
\begin{Notation} Let $ \gamma_i= \BBE (X_0X_i) $ and  $ \gamma^{(2)}_i= \BBE (X^2_0X_i) $. Define 
$\beta_{2, {\ell_k}} = 2 \sum_{i =1}^{{\ell_k}} \gamma_i$, $ \beta_{2}^{({\ell_k})} = 2 \sum_{i \geq {\ell_k} +1}  \gamma_i$ and   
$\beta_{3, 1,  {\ell_k}}= 3 \sum_{i = 1 }^{{\ell_k}} \ \gamma^{(2)}_i$. 
\end{Notation}

Next note that, since  $\BBE( X_{k-i}X_k ) = \gamma_i$, 
\begin{multline} \label{dec2bisbis}     
\frac{1}{2}  \BBE \{  f''_{n-k} (\S_{k-1})  \}  \beta_{2, {\ell_k}}  - A_{k,2}^{(1)}  
=    \sum_{i=1}^{\ell_k}  \gamma_i  \sum_{j=1}^i \BBE \{  f''_{n-k} (\S_{k-j})  -  f''_{n-k}  (\S_{k-j -1})   \}     
=  \\   \sum_{i=1}^{\ell_k}  \gamma_i \sum_{j=1}^i \BBE  \big \{   f^{(3)}_{n-k}  (\S_{k-j -1})   \X_{k-j}  \big \}     + 
  \sum_{i=1}^{\ell_k}  \frac{\gamma_i }{2}  \sum_{j=1}^i \BBE  \big \{   f^{(4)}_{n-k}  (\S_{k-j -1})   \X^2_{k-j}  \big \}   +  r_{n,k,2}^{(1)}  \, , \end{multline}
where 
\[
 r_{n,k,2}^{(1)} :=  \frac 12    \int_0^1 (1-t)^2  \sum_{i=1}^{\ell_k}  \gamma_i  \sum_{j=1}^i \BBE  \Big \{    f^{(5)}_{n-k}  (S_{k-j -1} + t \X_{k-j}) \X_{k-j}^3  \Big \}    dt  \, .
 \]
By Item 1 of Lemma  \ref{lmacrucial}, it follows that 
 \[
| r_{n,k,2}^{(1)} | \ll   M^4   (n-k+1)^{-3/2}   \sum_{i=1}^{\ell_k}  i  \theta(i )  \, .
 \]
Since $\sum_{i \geq 1}  i  \theta(i)   < \infty$, this implies that 
\beq \label{Boundr21}  
\sum_{k\in [1,n]} \vert  r_{n,k,2}^{(1)}  \vert   \ll 1 \, .
\eeq
Next, taking into account  Item 1 of Lemma  \ref{lmacrucial}, we get 
\[
  \sum_{i=1}^{\ell_k}  \sum_{j=1}^i \Big |  \BBE  \big \{   f^{(3)}_{n-k}  (\S_{k-j -1})   \BBE_0 (X_{k-j}  )  \big \}  \gamma_i  \Big |  \ll (n-k+1)^{-1/2}  \theta( k - j)  \theta(i ) \, .
\]
Hence, since $\ell_k = [k/2]$ and $\sum_{i \geq 1}  i  \theta(i)   < \infty$,
\begin{equation} \label{additional-20-1}   
\sum_{k=1}^n \sum_{i=1}^{\ell_k}  \sum_{j=1}^i \Big |  \BBE  \big \{   f^{(3)}_{n-k}  (\S_{k-j -1})   \BBE_0 (X_{k-j}  )  \big \}  \gamma_i  \Big |  \ll  \sum_{k=1}^n   \theta( [k/2])   \sum_{i=1}^{ \ell_k }  i  \theta(i )   \ll 1  \, .
\end{equation}
With similar arguments, we have
\begin{equation} \label{additional-20-2}   
\sum_{k=1}^n \sum_{i=1}^{\ell_k}  \sum_{j=1}^i \Big |  \BBE  \big \{   f^{(4)}_{n-k}  (\S_{k-j -1})   (\X^2_{k-j}   - X^2_{k-j}  )  \big \}  \gamma_i  \Big |    \ll 1  \, .
\end{equation}
In addition, by taking into account  Item 2 of Lemma  \ref{lmacrucial}, we get 
\[
\Big |  \BBE  \big \{   f^{(4)}_{n-k}  (\S_{k-j -1})   \BBE (X^2_{k-j}  )  \big \}  \gamma_i  \Big |  \ll M^3  ( (n-k+1)^{-3/2} + (n-j)^{-1})   \theta (i) \, .
\]
Hence,
\begin{multline}  \label{additional-20-3} 
\sum_{k=1}^n\sum_{i=1}^{\ell_k}   \sum_{j=1}^i  \Big |  \BBE  \big \{   f^{(4)}_{n-k}  (\S_{k-j -1})   \BBE (X^2_{k-j}  )  \big \}  \gamma_i  \Big |  \\
  \ll \sum_{k=1}^n ( (n-k+1)^{-3/2} + (n-{\ell_k} )^{-1})    \sum_{i=1}^{\ell_k}  i \ \theta(i)  \ll \sum_{i=1}^n i \theta(i) \ll 1  \, .
 \end{multline}
So overall, starting from \eqref{dec2bisbis}  and taking into account the upper bounds \eqref{Boundr21}-\eqref{additional-20-3}, we derive that 
\begin{multline} \label{dec2bisbis-20sept}     
\frac{1}{2}  \BBE \{  f''_{n-k} (\S_{k-1})  \}  \beta_{2, {\ell_k}}  - A_{k,2}^{(1)}  =
    \sum_{i=1}^{\ell_k}  \gamma_i   \sum_{j=1}^i \BBE  \big \{   f^{(3)}_{n-k}  (\S_{k-j -1})   X_{k-j}  \big \}    \\ + 
  \sum_{i=1}^{\ell_k}  \frac{\gamma_i}{2} \sum_{j=1}^i \BBE  \big \{   f^{(4)}_{n-k}  (\S_{k-j -1})   (X^2_{k-j}  )^{(0)}\big \}  +  R_{n,k,2}^{(1)}  \, , \end{multline}
where $ R_{n,k,2}^{(1)} $ is such that 
\beq \label{BoundR21-20sept}  
\sum_{k\in [1,n]} \vert  R_{n,k,2}^{(1)}  \vert   \ll 1 \, .
\eeq

Now, let $ r_{n,k,3}^{(1)} = \frac{1}{3}  \BBE \{  f^{(3)}_{n-k} (\S_{k-1})  \}  \beta_{3,1, {\ell_k}}  - A_{k,3}^{(1)}$.
Then, recalling the notation $\gamma_i^{(2)} =   \BBE( X_0^2 X_i )$,  
\begin{multline} \label{dec2bisbisavec3}   
 r_{n,k,3}^{(1)}= \sum_{i=1}^{{\ell_k}}   \gamma_i^{(2)} \BBE \{  f^{(3)}_{n-k} (\S_{k-1})  \}   - A_{k,3}^{(1)} 
=    \sum_{i=1}^{\ell_k}  \gamma_i^{(2)} \sum_{j=1}^i \BBE  \big \{ f^{(3)}_{n-k} (\S_{k-j})  -   f^{(3)}_{n-k}  (\S_{k-j -1})  \big   \}   \\
=    \sum_{i=1}^{\ell_k}  \gamma_i^{(2)}  \sum_{j=1}^i  \Big (  \BBE  \big \{    f^{(4)}_{n-k}  (\S_{k-j -1} ) \X_{k-j}  \big \}     +   \int_0^1  (1-t)  \BBE  \big \{    f^{(5)}_{n-k}  (\S_{k-j -1} + t \X_{k-j}) \X^2_{k-j}  \big \}   dt  \Big )  \\  :=  r_{n,k,3}^{(1)}  (1) +  r_{n,k,3}^{(1)} (2)  \, .
\end{multline}
Taking into account Item 1 of Lemma \ref{lmacrucial} and the fact that $| \gamma_i^{(2)}  | \leq M^2   \theta(i)$ and $ \Vert   \X_{k-j} \Vert_{\infty}  \leq 2M$, it follows that 
\[
\vert  r_{n,k,3}^{(1)} (2)  \vert  \ll  \sum_{i=1}^{\ell_k}   \sum_{j=1}^i 
 \frac{  \theta (i)  } {  (n-k+1)^{3/2} } 
 \ll    
 \frac{  \sum_{i=1}^{\ell_k}  i \ \theta (i) } { (n-k+1)^{3/2} } \, . 
 \]
Therefore, since $\sum_{i \geq 1} i  \theta(i)   < \infty$, 
\beq \label{Boundr312}    
\sum_{k\in [1,n]} \vert  r_{n,k,3}^{(1)}  (2) \vert   \ll 1 \, .
\eeq
On another hand,  by Item 1 of Lemma \ref{lmacrucial}, 
\[
 \big |  \gamma_i^{(2)}  \BBE  \big \{    f^{(4)}_{n-k}  (\S_{k-j -1} ) \BBE_0(X_{k-j}  )  \big \}  \big |  \ll  (n-k+1)^{-1} \theta(k-j)   \theta(i)  \, .
\]
Hence, since $\sum_{i \geq 1} i  \theta(i)   < \infty$,
\beq \label{new-reste132-20sept}  
 \sum_{k=1}^n \sum_{i=1}^{\ell_k}  \sum_{j=1}^i  \big |  \gamma_i^{(2)}  \BBE  \big \{    f^{(4)}_{n-k}  (\S_{k-j -1} ) \BBE_0(X_{k-j}  )  \big \}    \big |    \ll \sum_{k=1}^n \theta([k/2])  \sum_{i=1}^{n}  i \theta(i)  \ll  1 \, .
\eeq

Starting from  \eqref{firstdecomp}  and taking into account   \eqref{reste1}, \eqref{linddec1},   \eqref{tildeB}, \eqref{newtermorder4}, \eqref{linddec2},   \eqref{BoundB1}, \eqref{dec2bisbis-20sept}-\eqref{new-reste132-20sept}  and the fact that  $\beta_2 = \sigma^2 = \BBE (X_0^2) + \beta_{2, \ell_k} + \beta_{2}^{(\ell_k)}$, we get 
\begin{align} \label{linddec3}
\BBE (   \Delta_{n,k} )   
 & =  \BBE  (  \Delta_{k,2}^{(1)}  ) + { \frac{1}{2} }  \BBE \big (   f''_{n-k} (\S_{k-1} ) (X_k^2 )^{(0)}\big )  -  { \frac{1}{2} }  \BBE \big (   f''_{n-k} (\S_{k-1} ) \big )  \beta_2^{({\ell_k})} 
 +
    \frac{1}{2}   \BBE  (\Delta_{k,3}^{(1)} )  \nonumber  \\
  &  -  \sum_{i=1}^{\ell_k}  \gamma_i \sum_{j=1}^i \BBE  \Big \{   f^{(3)}_{n-k}  (\S_{k-j -1})   X_{k-j}  \big \}   
  -    \sum_{i=1}^{\ell_k}   \frac{\gamma_i}{2} \sum_{j=1}^i \BBE  \Big \{   f^{(4)}_{n-k}  (\S_{k-j -1})   (X^2_{k-j} )^{(0)} \big \}   \nonumber   \\ 
  & -   \sum_{i=1}^{\ell_k}   \frac{\gamma_i^{(2)}}{2}  \sum_{j=1}^i  \Big (  \BBE  \big \{    f^{(4)}_{n-k}  (\S_{k-j -1} ) X_{k-j}  \big \}     + 
\frac{1}{6}   \BBE \big (  f_{n-k}^{(3)} (\S_{k-1} )  (X_k^3 - (  \beta_3 -  \beta_{3,1, {\ell_k}}  ) ) \big )  \nonumber   \\
& +   \frac{1}{6}   \BBE  (\Delta_{k,4}^{(1)} ) + \Gamma^{(1)}_{n,k} \, ,
\end{align}
where $  \Gamma^{(1)}_{n,k}   $ satisfies
\beq \label{boundGamma}
\sum_{k\in [1,n]} \vert  \Gamma^{(1)}_{n,k}   \vert  \ll 1 \, .
\eeq
Note first that 
\beq\label{B0term1}
\sum_{k=1}^n   \vert \BBE \big (   f''_{n-k} (\S_{k-1} ) \big )  \beta_2^{({\ell_k})}  \vert  \ll \sum_{k=1}^n   \vert   \beta_2^{({\ell_k})}  \vert \ll \sum_{k=1}^n  \sum_{i \geq {\ell_k}  +1}  \theta (i)  \ll  
 \sum_{i \geq 1}  i  \theta(i)  \ll 1 \, .
\eeq
To  handle the first two terms in the right hand side of \eqref{linddec3},  define  
\beq  \label{defmki}
m_k = [\sqrt{n-k}], \ m_{k,i} = \min ( m_k , k-i-1) \ \text{and}\ D_{k,i,2}^{(1)} = \BBE   \big \{   f''_{n-k} (\S_{k-i- 1})    ( X_{k-i}X_k )^{(0)}   \big \} . 
\eeq 
Then, for any integer $i$ in $[0, {\ell_k} ]$,  with the convention that $\S_u =0$ for any $u \leq 0$,  we write
\[
D_{k,i,2}^{(1)}  =  \BBE   \Big \{   \Bigl( f''_{n-k} (\S_{k- i - m_{k,i} -1 })   +   \sum_{j=i+1}^{i+m_{k,i} } 
 (  f''_{n-k} (\S_{k-j }) -  f''_{n-k} (\S_{k-j -1})  \Bigr) 
( X_{k-i}X_k )^{(0)}   \Big \} . 
\]
Let then, for $\ell =3,4,5$ and $t$ in $[0,1]$, 
\beq \label{defdelta24i}
{\tilde  \Delta}_{k,i,2}^{(1,\ell)} (t)   = \sum_{j=i+1}^{i+m_{k,i} }        
 f^{(\ell)}_{n-k}(\S_{k-j -1}+ t\X_{k-j} )   \X_{k-j}^{\ell -2}  ( X_{k-i}X_k )^{(0)}     . 
\eeq
By the Taylor integral formula, 
\begin{align} \label{linddec4}
 D_{k,i,2}^{(1)}    = & \BBE   \Big \{   f''_{n-k} (\S_{k- i - m_{k,i} -1 })    ( X_{k-i}X_k )^{(0)}  
+  {\tilde  \Delta}_{k,i,2}^{(1,3)} (0) + \frac 12 {\tilde  \Delta}_{k,i,2}^{(1,4)} (0)     \Big \} \nonumber
\nonumber \\ &
 + \frac{1}{2}  
\int_0^1 (1-t)^2  \BBE   \big \{    {\tilde  \Delta}_{k,i,2}^{(1,5)} (t)   \big \} dt  .
\end{align}
But, since $\Vert    f''_{n-k} \Vert_{\infty} \ll 1$, 
\begin{multline*} 
\sum_{k=1}^n \sum_{i=0}^{{\ell_k}} \big  \vert \BBE   \big \{   f''_{n-k} (\S_{k- i - m_{k,i} -1 })    ( X_{k-i}X_k )^{(0)}   \big \}  \big  \vert  \ll  \sum_{k=1}^n \sum_{i=0}^{{\ell_k}}  
\big ( \theta (m_k)  + \theta (k-i)  \big )  \wedge \theta (i)  \\
\ll  \sum_{k=1}^n \Bigl( m_k  
\theta (m_k)  +  \sum_{i=m_k}^{\ell_k} \theta (i) +  k \theta([k/2]) \Bigr)
\ll  1+  \sum_{k=1}^{[\sqrt{n}]}  k^2  \theta(k) + n   \sum_{k \geq [\sqrt{n}]}   \theta(k) \, .
\end{multline*}
Hence
\begin{equation} \label{B1term1}
\sum_{k=1}^n \sum_{i=0}^{{\ell_k}} \big  \vert \BBE   \big \{   f''_{n-k} (\S_{k- i - m_{k,i} -1 })    ( X_{k-i}X_k )^{(0)}   \big \}  \big  \vert  \ll 1+  \sum_{k \geq 1}  k (k \wedge \sqrt{n}) \theta (k) \, .
\end{equation}
On another hand, by using Item 1 of Lemma \ref{lmacrucial}, 
\beq \label{B2term1}
\sum_{k=1}^n \sum_{i=0}^{{\ell_k}}  \big |  \BBE   \big \{      {\tilde  \Delta}_{k,i,2}^{(1,5)} (t)   \big \}  \big | 
  \ll  \sum_{k=1}^n \sum_{i=0}^{{\ell_k}}   \sum_{j=i+1}^{i+m_{k,i}} 
  \frac{\theta (j-i) \wedge \theta (i)}{(n-k+1)^{3/2}}
\ll \sum_{j=1}^nj \theta ([j/2])  \ll 1 \, .
\eeq
For $\ell=3,4$, set 
\beq \label{defdelta24i}
 \Delta_{k,i,2}^{(1,\ell)}   := \sum_{j=i+1}^{i+m_{k,i} }        f^{(\ell)}_{n-k}(\S_{k-j -1})   X_{k-j}^{\ell -2}  ( X_{k-i}X_k )^{(0)}     . 
\eeq
Applying  Item 1 of Lemma \ref{lmacrucial} and using that $m_{k,i} \leq \sqrt{n-k+1}$, we get 
\begin{multline} \label{linddec5ante1}
 \sum_{k=1}^n  \sum_{i=0}^{{\ell_k}}   \big |\BBE   \big \{      {\tilde  \Delta}_{k,i,2}^{(1,3)} (0)  -  \Delta_{k,i,2}^{(1,3)}  \big \} \big |
\leq  \sum_{k=1}^n  \sum_{i=0}^{{\ell_k}}   \sum_{j=i+1}^{i+m_{k,i} }  \Vert      \BBE_0 (X_{k-j}  )  f^{(3)}_{n-k}(\S_{k-j -1})   ( X_{k-i}X_k )^{(0)}   \Vert_1  \\
\ll    \sum_{k=1}^n  \sum_{i=0}^{{\ell_k}}   \sum_{j=i+1}^{i+m_{k,i} } 
 \frac{  \theta( k-j) \wedge  \theta( j-i) \wedge \theta (i) } {( n-k+1)^{1/2} }
  \ll  \sum_{k=1}^n k \theta ([k/3]) \ll 1 \, .
\end{multline}
Similarly, since $\Vert   f^{(4)}_{n-k}  \Vert_{\infty} \ll (n-k+1)^{-1}$, we derive
\begin{equation}\label{linddec5ante2}
 \sum_{k=1}^n  \sum_{i=0}^{{\ell_k}}   \big |\BBE   \big \{     {\tilde  \Delta}_{k,i,2}^{(1,4)} (0)  -  \Delta_{k,i,2}^{(1,4)}   \big \} \big |   \ll 1  \, .
\end{equation}
Starting from \eqref{linddec3} and taking into account \eqref{boundGamma},  \eqref{B0term1}, \eqref{linddec4}, \eqref{B1term1}, \eqref{B2term1}, \eqref{linddec5ante1} and \eqref{linddec5ante2} , we then derive that 
\begin{align} \label{linddec5}
\BBE (   \Delta_{n,k} )   
&  =  \frac 12 \sum_{i=0}^{\ell_k} (1 + {\bf 1}_{\{i \neq 0 \}}) \BBE   \Big \{   {  \Delta}_{k,i,2}^{(1,3)} + \frac 12 {  \Delta}_{k,i,2}^{(1,4)}     \Big \}    +
    \frac{1}{2}   \BBE  (\Delta_{k,3}^{(1)} )  \nonumber  \\
  &  -  \sum_{i=1}^{\ell_k}  \gamma_i \sum_{j=1}^i \BBE  \Big \{   f^{(3)}_{n-k}  (\S_{k-j -1})   X_{k-j}  \big \}   
  -    \sum_{i=1}^{\ell_k}   \frac{\gamma_i}{2} \sum_{j=1}^i \BBE  \Big \{   f^{(4)}_{n-k}  (\S_{k-j -1})   (X^2_{k-j} )^{(0)} \big \}   \nonumber   \\ 
  & -   \sum_{i=1}^{\ell_k}   \frac{\gamma_i^{(2)}}{2}  \sum_{j=1}^i  \Big (  \BBE  \big \{    f^{(4)}_{n-k}  (\S_{k-j -1} ) X_{k-j}  \big \}     + 
\frac{1}{6}   \BBE \big (  f_{n-k}^{(3)} (\S_{k-1} )  (X_k^3 - (  \beta_3 -  \beta_{3,1, {\ell_k}}  ) ) \big )  \nonumber   \\
& +   \frac{1}{6}   \BBE  (\Delta_{k,4}^{(1)} ) + \Gamma^{(2)}_{n,k} \, ,\end{align}
where $  \Gamma^{(2)}_{n,k}   $ satisfies $\sum_{k=1}^n \vert  \Gamma^{(2)}_{n,k}   \vert  \ll 1 +  \sum_{k \geq 1}  k (k \wedge \sqrt{n}) \theta (k) $.  Introduce now the following notations. 
%\[
% \beta_{3,1}^{({\ell_k}) } = 3 \sum_{i \geq  \ell_k+1 } \BBE ( X^2_0 X_i)$ 
% \]
\begin{Notation} Let $\beta_{3, 2,  m_k} = 3 \sum_{i = 1 }^{m_{k,0}} \BBE ( X_0 X^2_i)$,
 $\beta^*_{3, {\ell_k} , m_k}= 6 \sum_{i= 1}^{{\ell_k} } \sum_{j=  i+1}^{i + m_{k,i}} \BBE ( X_0 X_{j-i}X_{j})$.
 Next, let 
 ${\tilde \beta}_3^{({\ell_k},m_k)} = \beta_3  -  \beta_{3,1, {\ell_k}}- \big \{  \BBE (X_0^3) +  \beta_{3,2, m_k} + \beta^*_{3,{\ell_k},m_k}  \big \}$, 
where we recall that $m_k$ and $m_{k,i}$ have been defined in \eqref{defmki}.  
\end{Notation}
Since
\begin{align*}
{\tilde \beta}_3^{({\ell_k},m_k)} =  & \quad 3  \sum_{i > {\ell_k} }  \BBE ( X_0^2 X_i) 
 +  3 \sum_{i > m_{k,0} }   \BBE ( X_0 X^2_i)   \\
 & +  6 \sum_{i= 1}^{{\ell_k}} \sum_{j  > m_{k,i} } \BBE ( X_0 X_jX_{j+i}) 
+ 6 \sum_{i >{\ell_k} } \sum_{j \geq 1} \BBE ( X_0 X_jX_{j+i}) ,
\end{align*} 
by  Item 1 of Lemma \ref{lmacrucial}, 
\begin{multline*} 
\sum_{k =1}^n  \Vert    f^{(3)}_{n-k}  \Vert_{\infty} \big |  {\tilde \beta}_3^{({\ell_k},m_k)} \big |   \ll  \sum_{k =1}^n  (n-k+1)^{-1/2} \sum_{i \geq {\ell_k}  \wedge m_k }  \theta (i) \\
+  \sum_{k =1}^n  (n-k+1)^{-1/2} \Bigl(   \sum_{i= 1}^{{\ell_k}} \sum_{j \geq m_{k,i}+1}  \theta (j) \wedge   \theta (i) 
+ \sum_{i \geq {\ell_k}} \sum_{j \geq 1}  \theta (j) \wedge   \theta (i) \Bigr) \, .
\end{multline*}
By simple algebra, and since $\sum_{i \geq 1} i\theta(i)  < \infty$, we then derive that
\beq \label{Boundtildebeta}
\sum_{k =1}^n \Vert f^{(3)}_{n-k} \Vert_{\infty}  |{\tilde \beta}_3^{({\ell_k},m_k)} | \ll  1 +  \sum_{i \geq 1}  i ( i \wedge \sqrt{n})\theta(i)   \, .
\eeq
Next we shall first center the random variables $X_{k-j}  ( X_{k-i}X_k )^{(0)} $ 
appearing  in the quantity $ \Delta_{k,i,2}^{(1,3)} $. Using that  $\BBE    \big \{ X_{k-j}  ( X_{k-i}X_k )^{(0)}   \big \} = \BBE    \big \{ X_{k-j}   X_{k-i}X_k   \big \} $, an application of Item 2 of Lemma \ref{lmacrucial} gives 
\begin{multline}  \label{defJ1k}
J_{1,k}:= \Big \vert  \frac{1}{6} \BBE  \big \{     f^{(3)}_{n-k}(\S_{k-1})  \big \}  \beta^*_{3,{\ell_k},m_k}   -  \sum_{i=1}^{\ell_k} \sum_{j=i+1}^{i+m_{k,i} }   \BBE  \big \{      f^{(3)}_{n-k}(\S_{k-j -1})   \big \}   \BBE    \big \{ X_{k-j}  ( X_{k-i}X_k )^{(0)}   \big \}   \Big \vert  \\
\ll   \sum_{i=1}^{\ell_k} \sum_{j=i+1}^{i+m_{k,i} }  \big \vert  \BBE  \big \{     f^{(3)}_{n-k}(\S_{k-1}) -    f^{(3)}_{n-k}(\S_{k-j -1})   \big \}    \big \vert  ( \theta(j-i) \wedge   \theta(i) )  \, .
\end{multline}
Let us handle the quantity $\BBE  \big \{     f^{(3)}_{n-k}(\S_{k-1}) -    f^{(3)}_{n-k}(\S_{k-j -1})   \big \} $. By Taylor integral formula,
\begin{multline*} 
 \BBE  \big \{     f^{(3)}_{n-k}(\S_{k-1}) -    f^{(3)}_{n-k}(\S_{k-j -1})   \big \}  = \sum_{\ell=1}^j   \BBE \big \{   f^{(4)}_{n-k}(\S_{k-\ell-1})  \X_{k- \ell} \big \}   \\
+ \int_0^1 (1-t) \sum_{\ell=1}^j  \BBE \big \{   f^{(5)}_{n-k}(\S_{k-\ell-1} + t \X_{k- \ell} )  \X^2_{k- \ell} \big \} dt  \, .
\end{multline*}
By using Item 1 of Lemma \ref{lmacrucial} and noticing that $\theta(j-i) \wedge   \theta(i)  \leq  \theta([j/2]) $, we get   
\begin{align}  \label{computation1-21}
&  \sum_{k =1}^n   \sum_{i=1}^{\ell_k} \sum_{j=i+1}^{i+m_{k,i} } \sum_{\ell=1}^j   \big |  \BBE \big \{   f^{(5)}_{n-k}(\S_{k-\ell-1} + t \X_{k- \ell} )  \X^2_{k- \ell} \big \}   ( \theta(j-i) \wedge   \theta(i) ) \big |   \nonumber   \\  & \ll  \sum_{k =1}^n  (n-k+1)^{-3/2}  \sum_{i=1}^{\ell_k} \sum_{j=i+1}^{i+m_{k,i} } j    \theta([j/2])  \nonumber  \\ &  \ll  \sum_{k =1}^n  (n-k+1)^{-3/2}   \Big \{  \sum_{j=1}^{2 [\sqrt n]} j^2    \theta(j )  +  m_k   \sum_{ j \geq [\sqrt{n}]} j  \theta(j )  \Big \}     \ll  \sum_{i \geq 1}  i ( i \wedge \sqrt{n})\theta(i)  \, .
\end{align}
Next, by Item 1 of Lemma \ref{lmacrucial} again, 
\begin{align}  \label{computation2-21}
  \sum_{k =1}^n   \sum_{i=1}^{\ell_k} &  \sum_{j=i+1}^{i+m_{k,i} } \sum_{\ell=1}^j   \big |   \BBE \big \{   f^{(4)}_{n-k}(\S_{k-\ell-1})  \BBE_0(X_{k- \ell} )  \big \} ( \theta(j-i) \wedge   \theta(i) ) \big |  \nonumber  \\ 
 &   \ll  \sum_{k =1}^n  (n-k+1)^{-1}  \sum_{i=1}^{\ell_k} \sum_{j=i+1}^{i+m_{k,i} }  \sum_{\ell=1}^j  \theta(k-\ell)  \theta([j/2])  \nonumber  \\ 
 &    \ll  \sum_{k =1}^{[n/2]}  (n-k+1)^{-1}  \sum_{i=1}^{\ell_k} \sum_{j=i+1}^{i+m_{k,i} } j  \theta([j/2]) \nonumber \\ & \quad \quad +   \sum_{k =[n/2] +1}^{n }  (n-k+1)^{-1}  \sum_{i=1}^{[k/2]} \sum_{j=i+1}^{i+m_k} \sum_{\ell=1}^{[k/2] + m_k}  \theta(k-\ell)   \theta([j/2])  \, .
\end{align}
With the  computations as given in \eqref{computation1-21}
and the fact that 
\begin{multline}   \label{computation3-21}
  \sum_{k =[n/2] +1}^{n }  (n-k+1)^{-1}  \sum_{i=1}^{[k/2]} \sum_{j=i+1}^{i+m_k} \sum_{\ell=1}^{[k/2] + m_k}  \theta(k-\ell)   \theta([j/2])   \\
      \ll     \sum_{k =[n/2] +1}^{n }  (n-k+1)^{-1}   k  \theta( [k/4] )  \sum_{i=1}^{[k/2]} \sum_{j=i+1}^{i+m_k}   \theta([j/2])   \ll  \sum_{k \geq 1 } k  \theta(k)   \sum_{i \geq 1}  \theta(i)   \, ,
\end{multline}
we derive, overall,  that 
\begin{equation} \label{computation4-21}
J_{1,k} \leq  \sum_{i=1}^{\ell_k} \sum_{j=i+1}^{i+m_{k,i} }    \sum_{\ell=1}^j \big \vert     \BBE \big \{   f^{(4)}_{n-k}(\S_{k-\ell-1})  X_{k- \ell} \big \}     \big \vert  \theta([j/2])  + \Gamma^{(3)}_{n,k} \, ,
\end{equation}
where $  \Gamma^{(3)}_{n,k}   $ satisfies $ 1+ \sum_{k=1}^n \vert  \Gamma^{(3)}_{n,k}   \vert  \ll  \sum_{i \geq 1}  i ( i \wedge \sqrt{n})\theta(i) $.  Next, for $m_{k, \ell}$ defined in \eqref{defmki}, write
\begin{multline*} 
 \BBE \big \{   f^{(4)}_{n-k}(\S_{k-\ell-1})  X_{k- \ell} \big \}  \\ =   \BBE \big \{   f^{(4)}_{n-k}(\S_{k-\ell-m_{k,\ell} -1})  X_{k- \ell} \big \}  +  \sum_{u=\ell+1}^{\ell + m_{k,\ell}} \BBE
  \big \{   (  f^{(4)}_{n-k}(\S_{k-u} ) -   f^{(4)}_{n-k}(\S_{k-u-1})  X_{k- \ell} \big \}  \, , 
\end{multline*}
implying, by using Item 1 of Lemma \ref{lmacrucial}, that 
\begin{multline*} 
 \big |    \BBE \big \{   f^{(4)}_{n-k}(\S_{k-\ell-1})  X_{k- \ell} \big \}  \big |  \\ \ll   (n-k+1)^{-1}   (   \theta( m_k ) + \theta ( k-\ell ) ) +     (n-k+1)^{-3/2}   \sum_{u=\ell+1}^{\ell + m_{k,\ell}}  \theta ( u-\ell )   \, .
\end{multline*}
Hence
\begin{multline*} 
 \sum_{k=1}^n  \sum_{i=1}^{\ell_k} \sum_{j=i+1}^{i+m_{k,i} }  \sum_{\ell =1}^j \big |    \BBE \big \{   f^{(4)}_{n-k}(\S_{k-\ell-1})  X_{k- \ell} \big \}  \big |   \theta([j/2])  \\
   \ll   \sum_{k =1}^n   \sum_{i=1}^{\ell_k} \sum_{j=i+1}^{i+m_{k,i} } (n-k+1)^{-1}       \theta([j/2])  \Big \{  j  \theta( m_k )      +    \sum_{\ell =1}^j   \theta ( k-\ell )  \Big \}    \\ +  \sum_{k =1}^n   \sum_{i=1}^{\ell_k} \sum_{j=i+1}^{i+m_{k,i} } (n-k+1)^{-3/2}    j    \theta([j/2])  \sum_{u=1}^{\sqrt n}  \theta ( u)   \, .
\end{multline*}
With the  computations as given in \eqref{computation1-21}-\eqref{computation3-21} together with the fact that 
\[
 \sum_{k =1}^n  (n-k+1)^{-1} m_k  \theta( m_k )  = \sum_{k=1}^n k^{-1/2} \theta( [\sqrt{k}] )   \ll \sum_{k=1}^{[\sqrt n]}   \theta( k)  \, , 
\]
it follows that 
\begin{equation} \label{computation5-21}
 \sum_{k=1}^n  \sum_{i=1}^{\ell_k} \sum_{j=i+1}^{i+m_{k,i} }  \sum_{\ell =1}^j \big |    \BBE \big \{   f^{(4)}_{n-k}(\S_{k-\ell-1})  X_{k- \ell} \big \}  \big |   \theta([j/2])   \ll 1+  \sum_{i \geq 1}  i ( i \wedge \sqrt{n})\theta(i)   \, .
\end{equation}
Therefore \eqref{computation4-21}  together with \eqref{computation5-21} imply 
\begin{multline} \label{B1dec5}
 \sum_{k=1}^n \Big \vert  \frac{1}{6} \BBE  \big \{     f^{(3)}_{n-k}(\S_{k-1})  \big \}  \beta^*_{3,{\ell_k},m_k}   -  \sum_{i=1}^{\ell_k} \sum_{j=i+1}^{i+m_{k,i} }   \BBE  \big \{      f^{(3)}_{n-k}(\S_{k-j -1})   \big \}   \BBE    \big \{ X_{k-j}  ( X_{k-i}X_k )^{(0)}   \big \}   \Big \vert \\
\ll 1+  \sum_{i \geq 1}  i ( i \wedge \sqrt{n})\theta(i)  \, .
\end{multline} 
With similar arguments, we infer that 
\begin{align} \label{B1dec5i=0}
 \sum_{k=1}^n \Big \vert  \frac{1}{3} \BBE  \big \{     f^{(3)}_{n-k}(\S_{k-1})  \big \}  \beta_{3,2,m_k}   -  \sum_{j=1}^{m_{k,0} }   \BBE  \big \{      f^{(3)}_{n-k}(\S_{k-j -1})   \big \}  & \BBE    \big \{ X_{k-j}  ( X^2_k )^{(0)}   \big \}   \Big \vert
\nonumber \\ & 
\ll 1+  \sum_{i \geq 1}  i ( i \wedge \sqrt{n})\theta(i)  \, .
\end{align}
Now, for any integer $i \in [0,n]$, let 
\[
{ \Delta}_{k,i,2}^{(1,3,0)}   := \sum_{j=i+1}^{i+m_{k,i} }   \big \{     f^{(3)}_{n-k}(\S_{k-j -1})   (  X_{k-j}  ( X_{k-i}X_k )^{(0)}  )^{(0)}  \big \}  \, .
\]
Starting from \eqref{linddec5} and taking into account  \eqref{Boundtildebeta}, \eqref{B1dec5} and \eqref{B1dec5i=0}, we then obtain  
\begin{align} \label{linddec6}
 & \BBE (   \Delta_{n,k} )   
 =  \frac 12 \sum_{i=0}^{\ell_k} (1 + {\bf 1}_{\{i \neq 0 \}}) \BBE   \big \{   {  \Delta}_{k,i,2}^{(1,3,0)} \big \}   + 
    \frac{1}{2}   \BBE  (\Delta_{k,3}^{(1)} )    + 
\frac{1}{6}   \BBE \big (  f_{n-k}^{(3)} (\S_{k-1} )  (X_k^3 - \BBE (X_0^3) ) \big ) \nonumber  \\
 &   -  \sum_{i=1}^{\ell_k}  \gamma_i \sum_{j=1}^i \BBE  \big \{   f^{(3)}_{n-k}  (\S_{k-j -1})   X_{k-j}  \big \}   
  -   \sum_{i=1}^{\ell_k} \frac{\gamma_i}{2}   \sum_{j=1}^i \BBE  \big \{   f^{(4)}_{n-k}  (\S_{k-j -1})   (X^2_{k-j} )^{(0)} \big \}  \nonumber   \\
 &   -   \sum_{i=1}^{\ell_k} \frac{\gamma_i^{(2)}}{2}   \sum_{j=1}^i \BBE  \big \{   f^{(4)}_{n-k}  (\S_{k-j -1})   X_{k-j} \big \}  +  \frac 14 \sum_{i=0}^{\ell_k} (1 + {\bf 1}_{\{i \neq 0 \}}) \BBE   \big \{ {  \Delta}_{k,i,2}^{(1,4)}     \big \} +    \frac{1}{6}   \BBE  (\Delta_{k,4}^{(1)} )   + \Gamma^{(4)}_{n,k} \, ,
\end{align}
where $  \Gamma^{(4)}_{n,k}   $ satisfies $\sum_{k=1}^n \vert  \Gamma^{(3)}_{n,k}   \vert  \ll 1 + \sum_{i \geq 1}  i ( i \wedge \sqrt{n})\theta(i) $.  

In what follows we continue the estimation of each term in the right-hand side of \eqref{linddec6} and show that the sum over $k $ from $1$ to $n$ of their absolute values is bounded by a constant times $ \{1 + \sum_{i \geq 1}  i ( i \wedge \sqrt{n})\theta(i) \}$. Let us start by dealing with the quantities $  {  \Delta}_{k,i,2}^{(1,3,0)}$.  With this aim, note first that for $m_{k,j}$ defined in \eqref{defmki},   
\begin{multline*}
 \big |\BBE \big \{     f^{(3)}_{n-k}(\S_{k-j -m_{k,j} -1 })   (  X_{k-j}  ( X_{k-i}X_k )^{(0)}  )^{(0)}  \big \}  \big |
 \ll   \Vert  f^{(3)}_{n-k} \Vert_{\infty}  \big ( \theta (m_{k,j}) \wedge  \theta (j-i)    \wedge  \theta (i)  \big ) \\
 \ll \Vert  f^{(3)}_{n-k} \Vert_{\infty}  \big (   \theta (m_{k})  \wedge  \theta (j-i)    \wedge  \theta (i)  + \theta (k-j)  \wedge  \theta (j-i)    \wedge  \theta (i)   \big )    \, .
\end{multline*}
Hence, by Item 1 of Lemma \ref{lmacrucial} and the fact that $m_k\leq \sqrt{n-k+1}$,
\begin{multline} \label{B1dec6}
\sum_{k=1}^n \sum_{i = 0}^{\ell_k}  \sum_{j=i+1}^{i+m_{k,i} }   \big | \BBE \big \{     f^{(3)}_{n-k}(\S_{k-j -m_{k,j} -1 })   (  X_{k-j}  ( X_{k-i}X_k )^{(0)}  )^{(0)}  \big \}   \big |\\
 \ll  \sum_{k=1}^n    \frac{m_k}{ \sqrt{n-k+1}} 
 \Bigl(  m_k \theta (m_k)   +   \sum_{i \geq m_k} \theta (i)  +  \ell_k  \theta ([k/3] \Bigr)  \\
  \ll  \sum_{k=1}^n    \sqrt{k}  \theta ( [\sqrt k])   +  \sum_{k=1}^n    \sum_{i \geq  [\sqrt k]} \theta (i)  +  \sum_{k=1}^n   k  \theta ( k ) \ll 1 +  \sum_{i \geq 1}  i ( i \wedge \sqrt{n})\theta(i)    \, .
\end{multline}
On another hand, by the Taylor integral formula, 
\begin{align*}
\BBE \big \{   & \big ( f^{(3)}_{n-k}(\S_{k-j -1 })  -     f^{(3)}_{n-k}(\S_{k-j -m_{k,j} -1 })   \big )  (  X_{k-j}  ( X_{k-i}X_k )^{(0)}  )^{(0)}  \big \}  \\
&= \sum_{u=1}^{m_{k,j}}  \BBE \big \{  \big ( f^{(3)}_{n-k}(\S_{k-j -u })  -     f^{(3)}_{n-k}(\S_{k-j -u -1 })   \big )  (  X_{k-j}  ( X_{k-i}X_k )^{(0)}  )^{(0)}  \big \}   \\
&=  \sum_{u=1}^{m_{k,j}}  \BBE \big \{   f^{(4)}_{n-k}(\S_{k-j -u -1 }    )   \X_{k-j -u }  (  X_{k-j}  ( X_{k-i}X_k )^{(0)}  )^{(0)}  \big \}    \\
& \quad + \sum_{u=1}^{m_{k,j}}  \int_0^1 ( 1-t) \BBE \big \{   f^{(5)}_{n-k}(\S_{k-j -u -1 }  +t  \X_{k-j -u }   )   \X^2_{k-j -u }  (  X_{k-j}  ( X_{k-i}X_k )^{(0)}  )^{(0)}  \big \} dt   \, .
\end{align*}
According to Item 1 of Lemma \ref{lmacrucial}, 
\begin{align} \label{B2dec6}
\sum_{k=1}^n \sum_{i = 0}^{\ell_k} & \sum_{j=i+1}^{i+m_{k,i} }  \sum_{u=1}^{m_{k,j}}  | \BBE \big \{   f^{(5)}_{n-k}(\S_{k-j -u -1 }  +t  \X_{k-j -u }   )   \X^2_{k-j -u }  (  X_{k-j}  ( X_{k-i}X_k )^{(0)}  )^{(0)}  \big \}  |  \nonumber  \\
& \ll      \sum_{k=1}^n\Vert  f^{(5)}_{n-k} \Vert_{\infty}   \sum_{i = 0}^{\ell_k}  \sum_{j=1}^{m_k }  \sum_{u=1}^{m_k}  \big ( \theta (u) \wedge  \theta (j)    \wedge  \theta (i)  \big )  \nonumber  \\
&  \ll      \sum_{k=1}^n  (n-k+1)^{-3/2}  \Big \{ \sum_{u =1}^{[\sqrt n]}  u^2   \theta (u) +  m_k^2  \sum_{i \geq  m_k+1}     \theta (i)  \Big  \} 
  \nonumber   \\ & \ll  \sum_{u =1}^{[\sqrt n]}  u^2 \theta (u)  +  \sum_{i \geq 1}   i  \theta (i)  \ll 1 +  \sum_{i \geq 1}  i ( i \wedge \sqrt{n})\theta(i)  \, .
\end{align}
Next, let 
$Z_{k,j,u,i}:= \X_{k-u }  (  X_{k-j}  ( X_{k-i}X_k )^{(0)}  )^{(0)}$ and 
$Z^{(0)}_{k,j,u,i}:=Z_{k,j,u,i}-\BBE ( Z_{k,j,u,i})$. Since
\[
 \vert \BBE ( Z_{k,j,u,i})  \vert  \ll 
 \bigl(  \theta (u-j ) + \theta (k-j) \bigr) \wedge    \theta (j-i) \wedge   \theta (i) \, , 
 \]
by Item 2 of Lemma \ref{lmacrucial}, 
\begin{align}  \label{B3dec6ante}
& \sum_{k=1}^n   \sum_{i=0}^{\ell_k}   \sum_{j=i+1}^{ i+m_{k,i}}  \sum_{u= j+1}^{j+m_{k,j}}   \big | \BBE \big \{   f^{(4)}_{n-k}(\S_{k -u -1 }    )  \big \}  \BBE  ( Z_{k,j,u,i} )   \big |   \nonumber  \\
&\ll \sum_{k=1}^n  \sum_{i=0}^{\ell_k}   \sum_{j=i+1}^{ i+m_{k,i}}   \sum_{u=j+1}^{j+m_{k,j}}  
\frac{ \bigl( \bigl(  \theta (u-j ) + \theta (k-j) \bigr) \wedge    \theta (j-i) \wedge   \theta (i) \bigr) }
{( n-u ) \wedge (n-k+1)^{3/2} }
\nonumber \\ &
\ll \sum_{k=1}^n   \big (  (n-k+1)^{-3/2} + n^{-1} \big )  \Bigl( m^2_k \sum_{i =m_k }^{{\ell_k}} \theta (i)  +  \sum_{u=1}^{m_k} u^2  \theta (u) + m_k^2    k \theta ([k/3] \Bigr) 
 \nonumber \\ & 
 \ll 1 +  \sum_{u \geq 1}  u ( u \wedge {\sqrt n} )   \theta (u)    \, .
\end{align}
On another hand,  for $m_{k,u}$ defined in \eqref{defmki}, 
\begin{multline*}
 \big |   \BBE \big \{   f^{(4)}_{n-k}(\S_{k -u -m_{k,u} -1 }    )   Z^{(0)}_{k,j,u,i}   \big \}  \big |  
\ll   \Vert  f^{(4)}_{n-k} \Vert_{\infty}   \Big \{   \big ( \theta (m_k) \wedge  \theta (u-j)   \wedge \theta (j-i)    \wedge  \theta (i)  \big )   \\ +   \big ( \theta (k-u) \wedge  \theta (u-j)   \wedge \theta (j-i)    \wedge  \theta (i)  \big )   \Big \}   \, .
\end{multline*}
Hence, using Item 1 of Lemma \ref{lmacrucial} and the fact that $m_k^2 \leq n-k+1$, 
\begin{align} \label{B3dec6}
& \sum_{k=1}^n \sum_{i = 0}^{\ell_k}  \sum_{j=i+1}^{i+m_{k,i} }  \sum_{u=j+1}^{j+m_{k,j}}  
 \big |   \BBE \big \{   f^{(4)}_{n-k}(\S_{k-u -m_{k,u} -1 }    )   Z^{(0)}_{k,j,u,i}   \big \}  \big | 
\nonumber  \\
& \ll   \sum_{k=1}^n    \frac{m^2_k}{ n-k+1 }  \Bigl( m_k  \theta (m_k)   +   \sum_{i \geq m_k} \theta (i)  +   k    \theta ([k/4]  )\Bigr)   \nonumber  \\
& \ll  \sum_{k=1}^n    \Bigl( \sqrt{k}  \theta ([\sqrt k ])   +     \sum_{i \geq [\sqrt k ]} \theta (i)  +    k    \theta (k)  \Bigr)  \ll   1 +  \sum_{u \geq 1}  u ( u \wedge {\sqrt n} )   \theta (u)    \, .
\end{align}
Next
\begin{align*} 
  \BBE \big \{   \big (  f^{(4)}_{n-k}(\S_{k -u -1 }    )  & - f^{(4)}_{n-k}(\S_{k-u -m_{k,u} -1 }    )  \big )   Z^{(0)}_{k,j,u,i}   \big \}   
  \\
 & =
  \sum_{v=u+1}^{u+m_{k,u}}  \BBE \big \{   \big (  f^{(4)}_{n-k}(\S_{k -v }    ) - f^{(4)}_{n-k}(\S_{k-v -1 }    )  \big )   Z^{(0)}_{k,j,u,i}   \big \}    \\
&=  \sum_{v=u+1}^{u+m_{k,u}}  \int_0^1   \BBE \big \{   \big (  f^{(5)}_{n-k}(\S_{k -v -1 }  + t \X_{k-v }     )  \X_{k -v }   Z^{(0)}_{k,j,u,i}   \big \}  dt   \, .
\end{align*}
Therefore,  by Item 1 of Lemma \ref{lmacrucial}, 
\begin{align}  \label{B4dec6}
& \sum_{k=1}^n \sum_{i = 0}^{\ell_k}  \sum_{j=i+1}^{i+m_{k,i} }  \sum_{u=j+1}^{j+m_{k,j}}  \big |  \BBE \big \{   \big (  f^{(4)}_{n-k}(\S_{k-u -1 }    ) - f^{(4)}_{n-k}(\S_{k-u -m_{k,u} -1 }    )  \big )   Z^{(0)}_{k,j,u,i}   \big \}   \big |  \nonumber \\
&\ll  \sum_{k=1}^n \sum_{i = 0}^{\ell_k}  \sum_{j=i+1}^{i+m_{k,i} }  \sum_{u=j+1}^{j+m_{k,j}}  \sum_{v=u+1}^{u+m_{k,u}}  \Vert  f^{(5)}_{n-k} \Vert_{\infty}    \big ( \theta (v-u) \wedge  \theta (u-j)   \wedge \theta (j-i)    \wedge  \theta (i)  \big ) \nonumber  \\
& \ll \sum_{k=1}^n  \frac{1}{(n-k+1)^{3/2}}  \sum_{\ell =1}^{m_k}  \ell^3  \theta ( \ell )  +  \sum_{k=1}^n  \frac{m_k^3}{(n-k+1)^{3/2}}  \sum_{i \geq m_k +1}  \theta (i)   \nonumber \\
&\ll \sum_{k=1}^n  \frac{1}{k^{3/2}}  \sum_{\ell =1}^{[{\sqrt k}]}  \ell^3  \theta (\ell )  +   \sum_{k=1}^n   \sum_{i \geq  {[\sqrt k]}}  \theta (i)  
 \ll  \sum_{u \geq 1}  u ( u \wedge {\sqrt n} )   \theta (u) \, .
\end{align}
Taking into account \eqref{B1dec6},  \eqref{B2dec6},  \eqref{B3dec6ante},  \eqref{B3dec6} and \eqref{B4dec6}, it follows that 
\beq \label{B5dec6}
\sum_{k=1}^n \sum_{i=0}^{{\ell_k}}  \big |  \BBE ( {  \Delta}_{k,i,2}^{(1,3,0)}  ) \big |  \ll 1 +  \sum_{u \geq 1}  u ( u \wedge {\sqrt n} )   \theta (u)  \, .
\eeq    
With similar (but even simpler) arguments, we infer that the sum over $k$ from $1$ to $n$ of the second and third  terms in the right-hand side of \eqref{linddec6} are also  bounded by a constant times $ \{1 + \sum_{u \geq 1}  u ( u \wedge {\sqrt n} )   \theta (u) \} $. More precisely,      
\beq \label{B6dec6}
\sum_{k=1}^n  \big \{    \big |   \BBE  (\Delta_{n,k,3}^{(1)} ) \big |  +  \big |   \BBE \big (  f_{n-k}^{(3)} (\S_{k-1} )  (X_k^3 - \BBE(X_0^3) ) \big )    \big |    \big \}  \ll 1 + \sum_{u \geq 1}  u ( u \wedge {\sqrt n} )   \theta (u)\, .
\eeq
We deal now with the fourth term of the right hand side of \eqref{linddec6}. With this aim, recalling the definition \eqref{defmki} of $m_{k,j}$, note that 
\[
\big | \BBE  \Big \{   f^{(3)}_{n-k}  (\S_{k-j -m_{k,j} - 1})   X_{k-j}  \big \}   \big | \ll \Vert  f^{(3)}_{n-k} \Vert_{\infty}  \theta (m_{k,j})    \ll \Vert  f^{(3)}_{n-k} \Vert_{\infty}   \big ( \theta (m_{k}) + \theta (k-j)  \big )  \, .
\]
Hence,  by Item 1 of Lemma \ref{lmacrucial} and recalling the notation $\gamma_i=  \BBE( X_{0}X_i )$, we get 
\begin{align}  \label{B7dec6} 
\sum_{k=1}^n \sum_{i=1}^{\ell_k}  \sum_{j=1}^i  & \big | \gamma_i  \BBE  \big \{   f^{(3)}_{n-k}  (\S_{k-j -m_{k,j}- 1})   X_{k-j}  \big \}   \big | \nonumber   \\
& \ll   \sum_{k=1}^n  (n-k+1)^{-1/2} \sum_{i=1}^{\ell_k}  \sum_{j=1}^i    \theta (i)   (  \theta (m_k)  +\theta(k-j) )  \nonumber  \\ & \ll     \sum_{k=1}^n  k^{-1/2} \theta ([\sqrt k])   \sum_{i\geq 1}i  \theta (i) +  \sum_{k=1}^n \theta([k/2])   \sum_{i \geq 1}  i  \theta (i)   \ll 1 \, .
\end{align}
Next, by the Taylor integral formula, 
\begin{multline*} 
\!\!\!\! \BBE  \big \{  \big (   f^{(3)}_{n-k}  (\S_{k-j -1})  -  f^{(3)}_{n-k}  (\S_{k-j -m_{k,j}- 1})   \big )  X_{k-j}  \big \}    
 =  \sum_{u=j+1}^{j+m_{k,j}}  \BBE  \big \{  \big (   f^{(3)}_{n-k}  (\S_{k-u})  -  f^{(3)}_{n-k}  (\S_{k-u- 1})   \big )  X_{k-j}  \big \} \\ \ \ 
=    \sum_{u=j+1}^{j+m_{k,j}} \Bigl(  \BBE  \big \{  f^{(4)}_{n-k}    (\S_{k -u- 1})  \X_{k-u} X_{k-j}  \big \}  
+    \int_0^1 (1-t)  \BBE  \big \{    f^{(5)}_{n-k}    (\S_{k -u- 1} + t  \X_{k-u})  \X^2_{k -u} X_{k-j}  \big \}  dt  \Bigr) \, .
\end{multline*}
But,  by Item 1 of Lemma \ref{lmacrucial}, 
\begin{multline}  \label{B8dec6} 
\sum_{k=1}^n \sum_{i=1}^{\ell_k}  \sum_{j=1}^i   \sum_{u=j+1}^{j+m_{k,i}} \big |  \gamma_i  \BBE  \big \{    f^{(5)}_{n-k}    (\S_{k -u- 1} + t  \X_{k -u})  \X^2_{k-u} X_{k-j}  \big \} \big |   \\
\ll    \sum_{k=1}^n  (n-k+1)^{-3/2}  \sum_{u \geq 1}  \theta (u)  \sum_{i \geq 1} i  \theta (i)   \ll 1  \, .
\end{multline}
On another hand,  by Item 1 of Lemma \ref{lmacrucial}  again, 
\begin{multline*} 
 \big |  \gamma_i   \BBE  \big \{  f^{(4)}_{n-k}    (\S_{k -u- 1}) ( \X_{k-u} - X_{k-u} )  X_{k-j}  \big \}   \big | =  \big |  \gamma_i   \BBE  \big \{  f^{(4)}_{n-k}    (\S_{k -u- 1}) \BBE_0(X_{k-u})  X_{k-j}  \big \}   \big | \\
 \ll (n-k+1)^{-1} ( \theta (k-u) \wedge \theta (u-j) ) \theta(i)  \ll (n-k+1)^{-1}  \theta ([(k-j)/2])  \theta(i)   \, .
\end{multline*}
Hence, 
\begin{align}   \label{B9dec6-cor-ante} 
\sum_{k=1}^n \sum_{i=1}^{\ell_k}  \sum_{j=1}^i   \sum_{u=j+1}^{j+m_{k,j}}  \big |  \gamma_i  &  \BBE  \big \{  f^{(4)}_{n-k}    (\S_{k -u- 1}) ( \X_{k-u} - X_{k-u} )  X_{k-j}  \big \}   \big |  \nonumber  \\
&  \ll \sum_{k=1}^n   \theta ([k/4]) \sum_{i=1}^n i  \theta (i)  \ll 1 \, .
\end{align}
Moreover, by  Item 2 of Lemma \ref{lmacrucial}, 
\begin{multline}  \label{B9dec6} 
\sum_{k=1}^n \sum_{i=1}^{\ell_k}  \sum_{j=1}^i   \sum_{u=j+1}^{j+m_{k,i}}   \big |  \gamma_i \BBE  \big \{  f^{(4)}_{n-k}    (\S_{k -u- 1})    \big \}  \BBE ( X_{k-u} X_{k-j} )  \big |  \\
 \ll  \sum_{k=1}^n \sum_{i=1}^{\ell_k}  \sum_{j=1}^i   \sum_{u=j+1}^{j+m_{k,i}} 
  \frac{ \theta (u-j)  \theta (i) } { (n-u) \wedge (n-k+1)^{3/2} }   
 \ll  \sum_{u \geq 1}  \theta ([u/2])  \sum_{i \geq 1} i  \theta (i)    < \infty \, .
\end{multline}
Hence, taking into account  \eqref{B7dec6}, \eqref{B8dec6}, \eqref{B9dec6-cor-ante} and  \eqref{B9dec6}, we derive that 
\begin{multline}  \label{B9bisdec6} 
\sum_{k=1}^n \sum_{i=1}^{\ell_k}  \sum_{j=1}^i  \big | \gamma_i  \BBE  \big \{   f^{(3)}_{n-k}  (\S_{k-j - 1})   X_{k-j}  \big \}    \big |   \\
\ll  1 +  \sum_{k=1}^n \sum_{i=1}^{\ell_k}  \sum_{j=1}^i    \sum_{u=j+1}^{j+m_{k,j}}  \big |  \gamma_i  \BBE  \big \{  f^{(4)}_{n-k}    (\S_{k-u- 1})  ( X_{k-u} X_{k-j} )^{(0)} \big \}   \big | \, .
\end{multline}
Next, recalling the definition \eqref{defmki} of $m_{k,u}$, by Item 1 of Lemma \ref{lmacrucial},  note that 
\begin{align}  \label{B10dec6} 
\sum_{k=1}^n  &\sum_{i=1}^{\ell_k}  \sum_{j=1}^i   \sum_{u=j+1}^{j+m_{k,j}}  \big |  \gamma_i \BBE  \big \{  f^{(4)}_{n-k}    (\S_{k -u- m_{k,u} - 1}) (X_{k -u} X_{k-j} )^{(0)}  \big \}  \big | \nonumber  \\
&\ll    \sum_{k=1}^n \sum_{i=1}^{\ell_k}  \sum_{j=1}^i   \sum_{u=j+1}^{j+m_{k,j}}  \Vert f^{(4)}_{n-k} \Vert_{\infty}  \big \{   \theta (m_{k,u})  \wedge \theta (u-j) \big \}    \theta (i) \nonumber   \\
& \ll    \sum_{k=1}^n \frac{m_k \theta (m_k) }{n-k+1}   \sum_{i =1}^n  i  \theta (i)  +  \sum_{k=1}^n \sum_{i=1}^{\ell_k}  \sum_{j=1}^i   \sum_{u=j+1}^{k } \big \{   \theta (k-u)  \wedge \theta (u-j) \big \}  \theta (i) \nonumber \\
& \ll  \Big ( \sum_{i =1}^n i  \theta (i)   \Big)^2  \ll 1 \, .
\end{align}
On another hand, 
\begin{align*}  
 \big |  \BBE  \big \{  \big ( f^{(4)}_{n-k}    (\S_{k -u- 1})  & -  f^{(4)}_{n-k}    (\S_{k-u- m_{k,u} - 1})  \big ) (X_{k -u} X_{k-j} )^{(0)}  \big \}  \big |   \\
& \leq   \sum_{v =u+1}^{u+m_{k,u}}  \big |  \BBE  \big \{  \big ( f^{(4)}_{n-k}    (\S_{k-v})  -  f^{(4)}_{n-k}    (\S_{k- v - 1})  \big ) (X_{k -u} X_{k-j} )^{(0)}  \big \}  \big |  \\
& \leq     \sum_{v =u+1}^{u+m_{k,u}}  \int_0^1 \big |  \BBE  \big \{  f^{(5)}_{n-k}    (\S_{k- v - 1} + t X_{k-v})    \X_{k -v} (X_{k -u} X_{k-j} )^{(0)}  \big \}  \big |  dt \, .
\end{align*}
Hence,  
\begin{align}  \label{B11dec6} 
 \big |  \gamma_i  \BBE  \big \{  \big ( f^{(4)}_{n-k}    (\S_{k -u- 1})  & -  f^{(4)}_{n-k}    (\S_{k-u- m_{k,u} - 1})  \big ) (X_{k -u} X_{k-j} )^{(0)}  \big \}    \big |   \nonumber  \\
& \ll      \Vert f^{(5)}_{n-k} \Vert_{\infty}  \sum_{v =u+ 1}^{u+ m_k}  \big (  \theta (v-u) \wedge  \theta (u-j) ) \theta (i) 
\, .
\end{align}
Taking into account \eqref{B10dec6} and \eqref{B11dec6} together with Item 1 of Lemma \ref{lmacrucial}, it follows that
\begin{align}  \label{B12dec6} 
\sum_{k=1}^n \sum_{i=1}^{\ell_k} &  \sum_{j=1}^i   \sum_{u=j+1}^{j+m_{k,j}}  \big |  \gamma_i  \BBE  \big \{  f^{(4)}_{n-k}    (\S_{k-u- 1})  (X_{k -u} X_{k-j} )^{(0)}  \big \}   \big |  \nonumber  \\
&  \ll   1+    \sum_{k=1}^n  (n-k+1)^{-3/2}   \sum_{v=1}^{m_k}  v   \theta (v)   \sum_{i =1}^n  i  \theta (i)  \ll 1  \, . 
\end{align}   
Starting from \eqref{B9bisdec6}  and taking into account the upper bound \eqref{B12dec6}, we get that  the sum over $k$ from $1$ to $n$ of the fourth term in the right-hand side of \eqref{linddec6} is  uniformly bounded as a function of $n$. More precisely,  
\beq  \label{B9bisdec6term5} 
\sum_{k=1}^n \sum_{i=1}^{\ell_k}  \sum_{j=1}^i  \big |  \gamma_i \BBE  \big \{   f^{(3)}_{n-k}  (\S_{k-j - 1})   X_{k-j}  \big \}     \big |   
\ll  1 \, .
\eeq    
Similar computations (even simpler since we deal with the fourth derivative rather than the third one) give the following upper bound concerning the quantities involved in the fifth and sixth terms
in the right-hand side of \eqref{linddec6}: 
\beq  \label{B9bisdec6term5-new} 
\sum_{k=1}^n \sum_{i=1}^{\ell_k}  \sum_{j=1}^i   \Big (   \big | \gamma_i \BBE  \big \{    f^{(4)}_{n-k}  (\S_{k-j -1})   (X^2_{k-j} )^{(0)} \big \}       \big |  +  \big | \gamma_i^{(2)} \BBE  \big \{    f^{(4)}_{n-k}  (\S_{k-j -1})   X_{k-j}  \big \}       \big |    \Big ) 
\ll  1 \, .
\eeq    
We deal now with the last terms in the decomposition  \eqref{linddec6} and show that 
\beq  \label{lastaimdec6} 
\sum_{k=1}^n  \sum_{i=0}^{\ell_k}   \big |   \BBE ( \Delta_{k,i,2}^{(1,4)}  ) \big |  \ll 1 \mbox{ and }  \sum_{k=1}^n \big | 
  \BBE  (\Delta_{k,4}^{(1)} )  \big | 
\ll  1 \, ,
\eeq
where we recall that $\Delta_{k,i,2}^{(1,4)}  $ and $\Delta_{k,4}^{(1)}$  have been respectively defined in   \eqref{defdelta24i} and \eqref{defdeltank4}.
With this aim, note first that, by Item 2 of Lemma \ref{lmacrucial}, 
\begin{align}  \label{lastaimdec6P1} 
\sum_{k=1}^n &  \sum_{i=0}^{\ell_k}   \sum_{j=i+1}^{ i+m_{k,i}}  \big | \BBE  \big \{   f^{(4)}_{n-k}  (\S_{k-j - 1})   \BBE  \big ( X^2_{k-j}  (X_{k-i} X_k)^{(0)} \big ) \big \} \big |  \nonumber \\
&\ll \sum_{k=1}^n  \sum_{i=0}^{\ell_k}   \sum_{j=i+1}^{ i+m_{k,i}}   \big (  (n-k+1)^{-3/2} + (n-j)^{-1} \big )  \big ( \theta (j-i) \wedge   \theta (i) \ \big ) \nonumber \\
&\ll \sum_{k=1}^n   \big (  (n-k+1)^{-3/2} + n^{-1} \big )  \Big ( m_k \sum_{i \geq m_k } \theta (i)  +  \sum_{u=1}^{m_k} u  \theta (u) \Big ) % \nonumber \\ &
 \ll \sum_{u \geq 1}  u  \theta (u)   \, .
\end{align}  
Next, let $W_{k,i,j} =  \big ( X^2_{k-j}  (X_{k-i} X_k)^{(0)} \big )^{(0)}$.
We start by noticing that 
\begin{multline*} 
\sum_{k=1}^n  \sum_{i=0}^{\ell_k}   \sum_{j=i+1}^{ i+m_{k,i}}  \big | \BBE  \big \{   \big ( f^{(4)}_{n-k}  (\S_{k-j -m_{k,j}- 1})    W_{k,i,j}  \big \} \big |\\
\ll   \sum_{k=1}^n  \sum_{i=0}^{\ell_k}   \sum_{j=i+1}^{ i+m_{k,i}}  \Vert f^{(4)}_{n-k} \Vert_{\infty} \big ( \theta (m_{k,j}) \wedge  \theta (j-i) \wedge   \theta (i) \ \big )  \, .
\end{multline*}  
But, $ \theta (m_{k,j})  =  \theta (m_{k})  \vee  \theta (k-j) $. Hence, using Item 1 of Lemma \ref{lmacrucial},
\begin{align}  \label{lastaimdec6P2} 
& \sum_{k=1}^n  \sum_{i=0}^{\ell_k}   \sum_{j=i+1}^{ i+m_{k,i}}  \big | \BBE  \big \{   \big ( f^{(4)}_{n-k}  (\S_{k-j -m_{k,j}- 1})    W_{k,i,j}  \big \} \big | \nonumber \\
& \ll   \sum_{k=1}^n  (n-k+1)^{-1}  \Big (m_k^2  \theta (m_k) + m_k  \sum_{i \geq m_k} \theta (i)  +
m_k  \,  k \theta([k/3])\Big )  \nonumber  \\  
 & \ll  \sum_{k=1}^n   \theta (m_k) + \sum_{k=1}^n   (n-k+1)^{-1/2}   \sum_{i \geq m_k} \theta (i)  + \sum_{k=1}^n k \theta(k)  
\ll   \sum_{u \geq 1}  u  \theta (u)
 \, .
\end{align}  
Next,  by Item 1 of Lemma \ref{lmacrucial}, we derive
\begin{align}  \label{lastaimdec6P3} 
\sum_{k=1}^n &  \sum_{i=0}^{\ell_k}   \sum_{j=i+1}^{ i+m_{k,i}}  \big | \BBE  \big \{  \big (  f^{(4)}_{n-k}  (\S_{k-j -1})  -  f^{(4)}_{n-k}  (\S_{k-j -m_{k,j}- 1})  \big )   W_{k,i,j}  \big \} \big | 
 \nonumber   \\
& \leq 
 \sum_{k=1}^n  \sum_{i=0}^{\ell_k}   \sum_{j=i+1}^{ i+m_k}  \sum_{u=j+1}^{j+m_{k,j}}  \int_0^1 \big | \BBE  \big \{  \big (  f^{(5)}_{n-k}  (\S_{k-u-1} + t  \X_{k-u})  \big )   \X_{k -u} W_{k,i,j}  \big \} \big |  dt  \nonumber  \\
& \ll    \sum_{k=1}^n  \sum_{i=0}^{\ell_k}   \sum_{j=i+1}^{ i+m_k}  \sum_{u=j+1}^{j+m_{k,j}}  \Vert f^{(5)}_{n-k} \Vert_{\infty} \big ( \theta (u-j) \wedge  \theta (j-i) \wedge   \theta (i) \ \big )  \nonumber  \\
& \ll  \sum_{k=1}^n  (n-k+1)^{-3/2}  \Big (  \sum_{u=1}^{m_k}  u^2  \theta (u) + m^2_k \sum_{i \geq m_k} \theta (i)  \Big ) 
\ll   \sum_{u \geq 1}  u  \theta (u)  \, .
\end{align}  
Putting together \eqref{lastaimdec6P1}, \eqref{lastaimdec6P2} and \eqref{lastaimdec6P3}, the first part of \eqref{lastaimdec6} follows. Similar (but simpler) arguments lead to the second part of \eqref{lastaimdec6}. Finally, starting from \eqref{linddec6} and taking into account the upper bounds \eqref{B5dec6}, \eqref{B6dec6}, \eqref{B9bisdec6term5}, \eqref{B9bisdec6term5-new}  and \eqref{lastaimdec6}, it follows that  $ \sum_{k=1}^n \big |  {\mathbb E} ( \Delta_{n,k} )  \big | \ll 1 + \sum_{k\geq 1} k ( k \wedge \sqrt{n} )  \theta (k)$, which combined with \eqref{sumdelta} implies \eqref{aim1thmain} and then proves Item (b) of the theorem.

\subsection{Proof of Lemma \ref{lmacrucial}}

Item 1 comes from the smoothing lemma 6.1 in \cite{DMR09}. To prove Item 2, we write
\[
 \big | \BBE (  f_{n-k}^{(i) } ({\tilde S}_{\ell -1}  ) ) - \BBE (  f_{n-k}^{(i) } ({ S}_{\ell -1} )   )  \big |  \leq   \Vert f_{n-k}^{(i+1)}  \Vert_{\infty}  \Vert  \BBE_0 ( S_{\ell -1} ) \Vert_1 \, .
\]
Hence, since $ \Vert  \BBE_0 ( S_{\ell -1} ) \Vert_1 \leq \sum_{k=1}^{\ell} \theta_{X,1,1} \ll 1$, using Item 1, we derive that for any positive integer $\ell$, 
\[
 \big | \BBE (  f_{n-k}^{(i) } ({\tilde S}_{\ell -1})   ) - \BBE (  f_{n-k}^{(i) } ({ S}_{\ell -1} )  )  \big | \ll (n-k+1)^{-(i-1)/2} \, .
 \]
Next, let $(G_i)_{i \geq 1}$ be a sequence of iid centered Gaussian random variables with variance $\sigma^2$ and independent of $(X_i,B_i,Z_i)_{i \geq 1}$ (recall that the random variables $(B_i)$ and $(Z_i) $ have been defined at the beginning of Section \ref{subproofTh}). Let $N_{k} 
= \sum_{i=1}^{k} G_i $. Write that 
\[
\BBE (  f_{n-k}^{(i) } (S_{\ell -1}  )  ) = \BBE ( f_{n-k}^{(i) } (S_{\ell -1} )   ) - \BBE (  f_{n-k}^{(i) } (N_{\ell -1}  )   ) +   \BBE (  f_{n-k}^{(i) } (N_{\ell -1} )  ) \, .
\]
Next, let $t_k = \sigma \sqrt{(n-k)/2+1}$ and let
$\varphi_{t^2_k}$ be the density of the law ${\mathcal N} (0,t_k^2)$. Denote also by $H_{k,n} =  \sum_{i=k+1}^n B_i$ and note that, by definition, $H_{k,n} $ is independent of $S_{\ell -1} $ and of $N_{\ell -1} $.  Note that 
\[
 \BBE ( f_{n-k}^{(i) } (N_{\ell -1}  )  ) =  \BBE (    f * \varphi_{t^2_k}^{(i) } (N_{\ell -1}  + H_{k,n} )   ) =   \BBE (   f * \varphi_{t^2_k + \sigma^2 (\ell-1)}^{(i) } ( H_{k,n} )   )  \, .
\]
Using Item 1, it follows that 
\[
  |  \BBE ( f_{n-k}^{(i) } (N_{\ell -1}  )  ) |   \ll (n-k+ \ell)^{- (i-2)/2} \, .
\]
On another hand
\begin{multline*}
 \BBE (  f_{n-k}^{(i) } (S_{\ell -1} )   ) - \BBE (  f_{n-k}^{(i) } (N_{\ell -1}  )    )  
=    \BBE (    f * \varphi_{t^2_k}^{(i) } (S_{\ell -1} +H_{k,n}  )   )   -  \BBE (    f * \varphi_{t^2_k}^{(i) } (N_{\ell -1}  + H_{k,n} )   )  \\
= \int_{\mathbb R} \BBE \big \{  f'  (S_{\ell -1} +H_{k,n} -u )  -  f'   (N_{\ell -1} + H_{k,n} -u ) \big \}  \varphi_{t^2_k}^{(i-1) }   (u) du   \, .
\end{multline*}
Since  $f \in \Lambda_2(E)$, $g:=f'$ is in $\Lambda_1(E)$ meaning that  $g:{\mathbb R}
\times  E \rightarrow {\mathbb R}$ is measurable wrt the $\sigma$-fields 
${\mathcal L} ( {\mathbb R} \times E) $ and ${\mathcal B} ({\mathbb R})$, 
 $g( \cdot, w) $ is $1$-Lipschitz  and $g(0,w)=0$  for any 
$w \in E$. Therefore, since it is assumed that $\sum_{k \geq 1} k \theta_{X,3,4}(k) < \infty$, one can use  Item a) of Theorem 3.1 in \cite{DR08} (see also  Theorem 1.1 in \cite{Pe05}) which entails that 
\[
\sup_{v \in {\mathbb R}} \big | \BBE ( f'  (S_{\ell -1}  + v ) ) - \BBE( f'  (N_{\ell -1}  +v  ) )    \big | \ll 1 \, .
\] 
Note that  Item a) of Theorem 3.1 in \cite{DR08} is stated for $g$ a Lipschitz function but following its proof one can show  that it holds also if $g$ belongs to $\Lambda_1(E)$.   
On another hand, $\varphi_{t^2_k}^{(i-1) }   (u)  = t_k^{-i} \varphi_{1}^{(i-1) }   (u/t_k)$.  Therefore
\[
 \big |  \BBE (  f_{n-k}^{(i) } (S_{\ell -1} )  ) - \BBE ( f_{n-k}^{(i) } (N_{\ell -1}  )   )   \big | \ll  t_k^{ 1-i}  \Vert \varphi_{1}^{(i-1) }  \Vert_1 \ll  t_k^{ 1-i}  \, .
\]
Putting together all the above upper bounds gives Item 2 of Lemma \ref{lmacrucial}.

\section{Annex: convergence of quantiles in the CLT} \label{Annex}

\setcounter{equation}{0}

In this section, we give an inequality involving the difference between the 
quantile of a normalized random variable and the quantile of a standard normal, and the Wasserstein distance of 
order $p$ between the corresponding laws. The main result of this section is Proposition \ref{PropVaR} below which is a key result to prove Corollary \ref{AppliVaR}.

\begin{Proposition} \label{PropVaR} Let $Z$ be a centered real-valued random variable satisfying
$\BBE ( Z^2) \leq 2$.
Let $F_Z$ denote the distribution function of $Z$ and $\Phi$ denote the distribution function of a standard normal $Y$. For any $p\geq 1$, let 
$$
K_p =  \int_0^1 \big| F_Z^{-1} (t) - \Phi^{-1} (t) \big|^p dt .
$$
Then, for any $u$ in $(0,1/2]$, 
$$
\big| F_Z^{-1} (1-u) - \Phi^{-1} (1-u) \big| \leq \max \Bigl( \Bigl( \frac{(p+1)e K_p }{uQ_{1,Y} (u)} \Bigr)^{1/(p+1)} , 
\Bigl( \frac{(p+1)e K_p }{u} \Bigr)^{1/p} \Bigr) . 
$$
\end{Proposition}
\noindent
{\bf Proof.}  Throughout the proof, $H_Y = 1 - \Phi$ and $Q_Y$ is the inverse function of $H_Y$. With these
notations, 
\beq \label{ExpressionQ1Y}
Q_{1,Y} (u) = u^{-1} \int_0^u Q_Y (t) dt = u^{-1} \BBE \bigl( Y \BBI_{Y\geq Q_Y (u) } \bigr) = 
\frac{ \exp ( - Q_Y^2 (u) /2  ) }{ \sqrt{2\pi} \, u } \, .
\eeq
We also set $H_Z = 1-F_Z$ and we denote by $Q_Z$ the generalized inverse function of $H_Z$.  From
\eqref{ExpressionQ1Y}, Proposition \ref{PropVaR} is equivalent to 
\beq \label{ProofVaR1}
| Q_Z (u) - Q_Y (u) | \leq \max \Bigl( \bigl( (p+1)e \sqrt{2\pi} \,  e^{Q_Y^2 (u) /2}   K_p  \bigr)^{1/(p+1)} , 
\bigl( (p+1)e K_p /u \bigr)^{1/p} \Bigr)  
\eeq
for $u\leq 1/2$. We start by proving \eqref{ProofVaR1} in the case $Q_Z (u) > Q_Y (u)$.
\par\smallskip\no
{\sl Proof of \eqref{ProofVaR1} in the case $Q_Z (u) > Q_Y (u)$. } 
 Let $\delta = Q_Z (u) - Q_Y (u)$ and let $\eta$ be the unique real
in $(0,u)$ such that  $Q_Y (u-\eta ) = Q_Y (u) + \delta = Q_Z (u)$. From the convexity of $Q_Y$ on $(0,1/2]$, 
\beq \label{ProofVaR2}
Q_Y (u-t\eta )  \leq Q_Y (u) + t \delta \ \text{ for any } t\in [0,1] .
\eeq
Moreover $Q_Z ( u - t \eta ) \geq Q_Z (u) \geq Q_Y (u) + \delta$ for $t$ in $[0,1]$, whence, using the change of
variables $s=u-t\eta$,  
\beq
\label{ProofVaR3}
K_p  \geq 
\int_{u-\eta}^u |Q_Z (s) - Q_Y (s) |^p ds \geq \eta \int_0^1 (\delta -\delta t)^p dt = \eta \delta^p / (p+1) .
\eeq
In view of the above inequality, we have to bound $\eta$ from below. In order to get a lower bound on $\eta$,
we will bound up $-Q'_Y$. From the definition of $Q_Y$, 
$$
-Q'_Y (s) = -1/H'_Y ( Q_Y (s) ) = \sqrt{2\pi} \exp ( Q_Y^2 (s) / 2 ) \leq \sqrt{2\pi} \exp ( (Q_Y (u) + \delta)^2 /2 ) 
$$
for any $s$ in $[u-\eta , u]$,  
\par\smallskip
We now separate two cases. If $\delta \leq \sqrt{2 + Q_Y^2 (u) } - Q_Y (u)$, 
$$
- Q'_Y (s) \leq \sqrt{2\pi} \exp ( (Q_Y (u) + \delta)^2 /2 ) \leq \sqrt{2\pi}  \exp (1 + Q_Y^2 (u) /2 ) 
$$
for any $s$ in $[u-\eta , u]$. Then
\beq
\label{Competadelta}
 Q_Y (u-\eta ) - Q_Y (u) \leq \eta e\sqrt{2\pi}  \exp ( Q_Y^2 (u) /2 ) . 
\eeq
In that case, putting the above lower bound on $\eta$ in \eqref{ProofVaR3}, we obtain that 
\beq
\label{ProofVaR4}
   \delta^{p+1} \leq (p+1)e \sqrt{2\pi} \,  e^{Q_Y^2 (u) /2} K_p .
\eeq
\par
If $\delta > \sqrt{2 + Q_Y^2 (u) } - Q_Y (u)$, let $\delta_0 =  \sqrt{2 + Q_Y^2 (u) } - Q_Y (u)$ and let $\eta_0$
be the real in $(0,u)$ such that $Q_Y (u-\eta_0) = Q_Y (u) + \delta_0$. Then 
$\eta \geq \eta_0$ and $(\delta_0, \eta_0)$ still satisfies \eqref{Competadelta}, from which 
\beq
\label{ProofVaR5}
\eta \geq \eta_0 \ \geq  (e\sqrt{2\pi})^{-1} \Bigl(  \sqrt{2 + Q_Y^2 (u) } - Q_Y (u) \Bigr)  \exp ( - Q_Y^2 (u) /2 )  .
\eeq
Putting this lower bound in \eqref{ProofVaR3}, we obtain that 
\beq
\label{ProofVaR6}
   \delta^p \leq (p+1)e  \, (K_p/ u)  \frac{ \sqrt{2\pi}\exp (Q_Y^2 (u) /2) u }{ \sqrt{2 + Q_Y^2 (u) } - Q_Y (u) } .
\eeq
Now, setting $u= H_Y (x)$,  
$$
\sup_{u\in (0,1/2]} \frac{ \sqrt{2\pi}\exp (Q_Y^2 (u) /2) u }{ \sqrt{2 + Q_Y^2 (u) } - Q_Y (u) } = 
\sup_{x\geq 0} \frac{ \sqrt{2\pi}\exp (x^2/2) H_Y (x) }{ \sqrt{2 + x^2 } - x } \leq 1 
$$
by an inequality on the Mills ratio of Komatu \cite{Ko55}. The two above inequalities imply that 
\beq
\label{ProofVaR7}
   \delta^p \leq (p+1)e  \, (K_p/ u) ,
\eeq
if $\delta \geq \sqrt{2 + Q_Y^2 (u)} - Q_Y (u)$. 
Combining \eqref{ProofVaR4} and \eqref{ProofVaR7}, we get \eqref{ProofVaR1} in the case $Q_Z (u)> Q_Y (u)$.
It remains to prove \eqref{ProofVaR1} in the case $Q_Z (u) < Q_Y (u)$. 
\par\smallskip\no
{\sl Proof of \eqref{ProofVaR1} in the case $Q_Z (u) < Q_Y (u)$. }  Let then $\delta = Q_Y (u) - Q_Z (u)$. 
From the assumptions 
$\BBE (Z)=0$, $\BBE (Z^2) \leq 2$ and the Tchebichef-Cantelli inequality, for any $x\leq 0$, 
$H_Z (x) \geq  x^2 / (2+x^2)$. This which implies that 
\beq \label{LowerBoundQZ}
Q_Z (u) \geq - \sqrt{2u/(1-u)}\ \text{ for any }  u\in (0,1) .
\eeq
In particular, for $u\leq 1/2$, $Q_Z(u) \geq - \sqrt{2} \geq - \sqrt{2 + Q_Y^2 (u)}$. Let then $\beta$ be the positive
real such that $Q_Y ( u+ \beta) = Q_Z (u)$. From \eqref{LowerBoundQZ}, 
$-Q'_Y (s) \leq \sqrt{2\pi} \exp ( 1 + Q_Y^2 (u) /2 )$ for any $s$ in $[u, u+\beta]$. It follows that 
\[
Q_Y (u+s) \geq Q_Y (u) - s \sqrt{2\pi} \exp ( 1 + Q_Y^2 (u) /2 )
\]
for any $s$ in $[0,\beta]$.  For $s=\beta$, the above inequality yields 
\[ 
\beta \geq (e\sqrt{2\pi})^{-1} \exp ( - Q_Y^2 (u)/2  ) \delta := \eta .
\]
With the above definition of $\eta$, for any $t$ in $[0,1]$, 
\[ 
Q_Y (u + t\eta ) \geq Q_Y (u) - t \delta \geq Q_Y (u) - \delta \geq Q_Z ( u+t\eta ) . 
\] 
Hence 
\[ Q_Y (u + t\eta ) - Q_Z (u + t\eta ) \geq ( 1-t) \delta  \]
for any $t$ in $[0,1]$. It follows that 
\[ 
K_p \geq \eta \int_0^1 | Q_Y (u + t \eta ) - Q_Z (u+t\eta) |^p dt \geq \eta \int_0^1 (1-t)^p \delta^p dt = 
\frac{ \eta \delta^p } { p+1 } .
\]
The  above inequality together  with the definition of $\eta$ then imply \eqref{ProofVaR4}, which completes the proof of \eqref{ProofVaR1}. 

\medskip

\noindent {\bf Proof of  Corollary \ref{AppliVaR}.} Recall that  from Item (b) of Theorem \ref{thW2} (see also Comment \ref{commentaftermainth}), under the assumptions of Corollary   \ref{AppliVaR}, $W_2 (P_{S_n / \sigma_n}, G_1)  =O(n^{-1/2})$. Hence, Item (a) comes from  an application of Proposition \ref{PropVaR} by taking into account  the fact that, if $Y$ is a standard normal r.v.,  there exists a positive constant $\eta$ such that 
\[
\inf_{u \in (0,1/2]} \frac{ Q_{1,Y} (u)}{ \sqrt{\ln (1/u)}}  \geq \eta \, .
\]
Indeed, $Q_{1,Y} (u) \sim_{u \rightarrow 0} \sqrt{2 \ln (1/u)}$ and $Q_{1,Y} (1/2) >0$. 

Item (b) follows again  from Item (b) of Theorem \ref{thW2}  together with Inequality (2.7) in \cite{Rio17VaR}.

 %%%%%%%%%%%%%%%%%%%%%%%%%%%%%%% BIBLIOGRAPHY %%%%%%%%%%%%%%%%%%%%%%

\end{document}